\documentclass[draft,reqno]{amsproc}
\usepackage{amssymb}
\usepackage{euscript}   %\EuScript{A}
\usepackage{mathrsfs} %\mathscr T
\usepackage{units}   %\nicefrac %kolejarski ulamek
\usepackage{color}

\makeatletter
\@namedef{subjclassname@2010}{%
\textup{2010} Mathematics Subject Classification}
\makeatother

\numberwithin{equation}{subsection}

\newtheorem{thm}{Theorem}[subsection]
\newtheorem{cor}[thm]{Corollary}
\newtheorem{lem}[thm]{Lemma}
\newtheorem{pro}[thm]{Proposition}

\newtheorem*{thm*}{Theorem}

\theoremstyle{remark}
\newtheorem{rem}[thm]{Remark}

\theoremstyle{definition}
\newtheorem{exa}[thm]{Example}
\newtheorem*{caut}{{\it Caution}}

\DeclareMathOperator{\D}{d}
\DeclareMathOperator{\dess}{{\mathsf{Des}}}
\DeclareMathOperator{\dzii}{{\mathsf{Chi}}}
\DeclareMathOperator{\E}{e}
\DeclareMathOperator{\koo}{{\mathsf{root}}}

\DeclareMathOperator{\paa}{{\mathsf{par}}}
\DeclareMathOperator{\sumo}{{\sum\nolimits^{\smalloplus}}}

\newcommand*{\borel}[1]{{\mathfrak B}(#1)}
\newcommand*{\cbb}{\mathbb C}
\newcommand*{\card}[1]{\mathrm{card}(#1)}
\newcommand*{\des}[1]{{\dess(#1)}}
\newcommand*{\dtf}[2]{{\EuScript D}_{#1,#2}^{\varphi}}
\newcommand*{\dz}[1]{{\EuScript D}(#1)}
\newcommand*{\dzi}[1]{\dzii(#1)}
\newcommand*{\dzin}[2]{\dzii^{\langle#1\rangle}(#2)}
\newcommand*{\dzn}[1]{{\EuScript D}^\infty(#1)}
\newcommand*{\ee}{\mathcal E}
\newcommand*{\escr}{{\mathscr{E}_V}}
\newcommand*{\gamkap}[1]{\{\gamma_n\}_{n=#1}^\infty}
\newcommand*{\gams}{\{\gamma_n\}_{n=0}^\infty}
\newcommand*{\gamsj}{\{\gamma_{n+1}\}_{n=0}^\infty}
\newcommand*{\gamsmj}{\{\gamma_{n-1}\}_{n=0}^\infty}
\newcommand*{\gamst}[1]{\{\gamma_n(#1)\}_{n=0}^\infty}
\newcommand*{\gamstj}[1]{\{\gamma_{n+1}(#1)\}_{n=0}^\infty}
\newcommand*{\Ge}{\geqslant}
\newcommand*{\Ham}{H}
\newcommand*{\hh}{\mathcal H}
\newcommand*{\I}{{\mathrm i}}
\newcommand*{\imc}{\mathfrak{I}_{_{C}}}
\newcommand*{\is}[2]{\langle#1,#2\rangle}
\newcommand*{\kk}{\mathcal K}

\newcommand*{\lambdab}{{\boldsymbol\lambda}}
\newcommand*{\Le}{\leqslant}
\newcommand*{\mm}{\mathscr M}
\newcommand*{\muf}{\mu_{\mathrm F}}
\newcommand*{\muk}{\mu_{\mathrm K}}
\newcommand*{\nbb}{\mathbb N}

\newcommand*{\ob}[1]{{\mathcal R}(#1)}
\newcommand*{\pa}[1]{\paa(#1)}
\newcommand*{\pcal}{{\EuScript P}}
\newcommand*{\rbb}{\mathbb R}
\newcommand*{\rec}{\mathfrak{R}_{_{C}}}
\newcommand*{\sac}{\mathscr{S}_{A, C}^{\varphi}}
\newcommand*{\slam}{S_{\boldsymbol \lambda}}
\newcommand*{\smalloplus}{\raise0pt\hbox{$\scriptscriptstyle \oplus$}}

\newcommand*{\Sti}{S}
\newcommand*{\supp}[1]{\mathrm{supp}(#1)}
\newcommand*{\tcal}{{\mathscr T}}

\newcommand*{\zbb}{\mathbb Z}

\begin{document}
   \title[A Non-hyponormal Operator Generating
Stieltjes Moment Sequences] {A Non-hyponormal Operator
Generating \\ Stieltjes Moment Sequences}
   \author[Z.\ J.\ Jab{\l}o\'nski]{Zenon Jan Jab{\l}o\'nski}
\address{Instytut Matematyki, Uniwersytet Jagiello\'nski,
ul.\ \L ojasiewicza 6, PL-30348 Kra\-k\'ow, Poland}
   \email{Zenon.Jablonski@im.uj.edu.pl}
   \author[I.\ B.\ Jung]{Il Bong Jung}
   \address{Department of Mathematics, Kyungpook National
University, Daegu 702-701, Korea}
   \email{ibjung@knu.ac.kr}
   \author[J.\ Stochel]{Jan Stochel}
\address{Instytut Matematyki, Uniwersytet Jagiello\'nski,
ul.\ \L ojasiewicza 6, PL-30348 Kra\-k\'ow, Poland}
   \email{Jan.Stochel@im.uj.edu.pl}
   \thanks{Research of the first
and the third authors was supported by the MNiSzW
(Ministry of Science and Higher Education) grant NN201
546438 (2010-2013). The second author was supported by
Basic Science Research Program through the National
Research Foundation of Korea (NRF) grant funded by the
Korea government (MEST) (2009-0093125).}
   \begin{abstract}
A linear operator $S$ in a complex Hilbert space $\hh$
for which the set $\dzn{S}$ of its $C^\infty$-vectors
is dense in $\hh$ and $\{\|S^n f\|^2\}_{n=0}^\infty$
is a Stieltjes moment sequence for every $f \in
\dzn{S}$ is said to generate Stieltjes moment
sequences. It is shown that there exists a closed
non-hyponormal operator $S$ which generates Stieltjes
moment sequences. What is more, $\dzn{S}$ is a core of
any power $S^n$ of $S$. This is established with the
help of a weighted shift on a directed tree with one
branching vertex. The main tool in the construction
comes from the theory of indeterminate Stieltjes
moment sequences. As a consequence, it is shown that
there exists a non-hyponormal composition operator in
an $L^2$-space (over a $\sigma$-finite measure space)
which is injective, paranormal and which generates
Stieltjes moment sequences. In contrast to the case of
abstract Hilbert space operators, composition
operators which are formally normal and which generate
Stieltjes moment sequences are always subnormal (in
fact normal). The independence assertion of Barry
Simon's theorem which parameterizes von Neumann
extensions of a closed real symmetric operator with
deficiency indices $(1,1)$ is shown to be false.
   \end{abstract}
   \maketitle
   \section{PRELIMINARIES}
   \subsection{Introduction} A linear operator
$S$ in a complex Hilbert space $\hh$ is said to {\em
generate} Stieltjes moment sequences if the set
$\dzn{S}$ of all its $C^\infty$-vectors is dense in
$\hh$ and $\{\|S^n f\|^2\}_{n=0}^\infty$ is a
Stieltjes moment sequence for every $f\in \dzn{S}$.
The celebrated Lambert characterization of
subnormality \cite{lam0} states that a (closed)
bounded linear operator is subnormal if and only if it
generates Stieltjes moment sequences. As shown in
\cite{StSz1,b-j-j-s}, this result remains true for
some classes of unbounded operators (see \cite{Con}
and \cite{StSz3,StSz1,StSz4,StSz2} for the foundations
of the theory of bounded and unbounded subnormal
operators). To the best of our knowledge, the only
known examples of non-subnormal operators generating
Stieltjes moment sequences are those coming from
formally normal ones\footnote{\;Formally normal
operators are always hyponormal but not necessarily
subnormal (see \cite{Cod,Sch,sto}).} (see
\cite[Section 3.2]{b-j-j-s} for a more detailed
discussion of this question). Unfortunately, the
operators so constructed, though closable, are not
closed. In the present paper we provide an example of
a non-hyponormal (and thus a non-subnormal) closed
paranormal operator $S$ which generates Stieltjes
moment sequences\footnote{\;Note that if $S$ is a
Hilbert space operator which generates Stieltjes
moment sequences, then the operator $S|_{\dzn{S}}$ is
paranormal; see \eqref{para}.} and which has the
property that $\dzn{S}$ is a core of any power $S^n$
of $S$ (see Example \ref{exnindeterm}). This is a
carefully constructed weighted shift on an enumerable
leafless directed tree (we refer the reader to
\cite{j-j-s} for the foundations of the theory of
weighted shifts on directed trees). As a byproduct, we
obtain an example of a paranormal operator which is
not hyponormal (see \cite{fur,b-j,j-j-s} for other
examples of this kind).

Using N-extremal measures (including the Friedrichs
one) of an indeterminate moment sequence as well as
some facts from moment theory which relate the
determinacy of sequences $\{a_n\}_{n=0}^\infty$ and
$\{a_{n+1}\}_{n=0}^\infty$, we construct a
non-hyponormal weighted shift on a directed tree
$\tcal_{\infty,\kappa}$ which generates Stieltjes
moment sequences (cf.\ Example \ref{exnindeterm}). The
$\tcal_{\infty,\kappa}$ is an enumerable leafless
directed tree which has only one branching vertex
denoted by $0$. If $\kappa < \infty$, then
$\tcal_{\infty,\kappa}$ has a root and $0$ belongs to
the $\kappa$th generation of the root; otherwise
$\tcal_{\infty,\kappa}$ is rootless. The weighted
shift so constructed does not satisfy the consistency
condition \eqref{con2} at $u=0$ and it has no
consistent system of measures (in the sense of
\cite{b-j-j-s}). The case of $\kappa=\infty$ is
especially interesting because it leads to an example
of a non-hyponormal composition operator in an
$L^2$-space over a $\sigma$-finite measure space which
generates Stieltjes moment sequences (cf.\ Theorem
\ref{compns}). In view of \cite{b-j-j-s2}, this
example is the first showing that Lambert's
characterization of subnormality of composition
operators (cf.\ \cite{lam1}) is no longer true in the
unbounded case. As proved in \cite{b-j-j-s2}, each
formally normal composition operator in an $L^2$-space
is normal. This means that an example of a
non-subnormal formally normal operator $N$ with dense
set of $C^\infty$-vectors $f$ having the property that
$\{\|N^n f\|^2\}_{n=0}^\infty$ is a Stieltjes moment
sequence, could not be realized as a composition
operator in an $L^2$-space.

Since our main example (Example \ref{exnindeterm})
depends heavily on some subtle properties of
indeterminate Stieltjes moment sequences, we provide
necessary facts concerning N-extremal measures
including Krein and Friedrichs ones (see Sections
\ref{sec3} and \ref{sec4}). In Section \ref{sec5} we
supply examples of exotic Stieltjes moment sequences
that are used in Example \ref{exnindeterm}. The
necessary facts concerning weighted shifts $\slam$ on
directed trees are given in Section \ref{sec6}. Powers
of such operators are described in Section \ref{sec7}.
As a consequence, it is shown that if $\dzn{\slam}$ is
dense in the underlying Hilbert space, then
$\dzn{\slam}$ is a core of any power $\slam^n$ of
$\slam$. A sufficient condition for $\slam$ to
generate Stieltjes moment sequences, written in terms
of basic vectors, is given in Theorem \ref{cinf}.
Section \ref{sec8} offers a general scheme for
constructing weighted shifts on the directed tree
$\tcal_{\eta,\kappa}$ with assorted properties (cf.\
Theorem \ref{import}). Section \ref{sec9} contains the
main example of the paper. The appendix shows that the
independence assertion of Barry Simon's theorem which
parameterizes von Neumann extensions of a closed real
symmetric operator with deficiency indices $(1,1)$ is
false (cf.\ Proposition \ref{BS5}). This theorem was
used by Simon to describe N-extremal measures of
indeterminate moment sequences in \cite{sim}.
Fortunately, this fault does not spoil\footnote{\;May
be with an exception of Remark 2 on page 104 in
\cite{sim}.} the main idea of his paper which is based
on the formula (4.20) in \cite{sim}.
   \subsection{Notation and terminology}
In what follows, $\cbb$, $\rbb$ and $\zbb$ stand for
the sets of complex numbers, real numbers and integer
numbers, respectively. Set
   \begin{align*}
\text{$\nbb=\{n \in \zbb\colon n \Ge 1\}$, $\zbb_+ =
\nbb \cup \{0\}$, $\rbb_+ = \{x \in \rbb \colon x \Ge
0\}$.}
   \end{align*}
For a Borel set $\varOmega$ in $\rbb_+$, we denote by
$\borel{\varOmega}$ the $\sigma$-algebra of all Borel
sets in $\varOmega$. Given $a\in \rbb_+$, we write
$\delta_a$ for the Borel probability measure on
$\rbb_+$ concentrated on $\{a\}$. The closed support
of a finite positive Borel measure $\mu$ on $\rbb$
will be denoted by $\supp{\mu}$. We write $\card{X}$
for the cardinal number of a set $X$.

Let $A$ be a (linear) operator in a complex Hilbert
space $\hh$. Denote by $\dz{A}$, $\ob{A}$, $\ker(A)$,
$\bar A$ and $A^*$ the domain, the range, the kernel,
the closure and the adjoint of $A$ (in case they
exist). Set $\dzn{A} = \bigcap_{n=0}^\infty \dz{A^n}$;
members of $\dzn{A}$ are called {\em
$C^\infty$-vectors}. A linear subspace $\ee$ of
$\dz{A}$ is said to be a {\em core} of $A$ if the
graph of $A$ is contained in the closure of the graph
of the restriction $A|_{\ee}$ of $A$ to $\ee$. We say
that $A$ is {\em symmetric} if $A$ is densely defined,
$\dz{A} \subseteq \dz{A^*}$ and $Af = A^*f$ for all
$f\in \dz{A}$. If $A$ is densely defined and $A=A^*$,
then $A$ is called {\em selfadjoint}. The operator $A$
is said to be {\em essentially} selfadjoint if $A$ is
closable and the closure of $A$ is selfadjoint. The
orthogonal dimensions of $\ker(A^* \mp \I I)$, which
are denoted by $d_{\pm}=d_{\pm}(A)$, are called the
{\em deficiency indices} of a symmetric operator $A$
($I$ is the identity operator on $\hh$). It is
well-known that if $A$ is symmetric, then $A$ is
essentially selfadjoint if and only if its deficiency
indices are both equal to $0$. If $A$ is symmetric,
then $A$ has equal deficiency indices if and only if
it has a selfadjoint extension in $\hh$; such an
extension will be called a {\em von Neumann} extension
of $A$. Note that a symmetric operator may have no von
Neumann extension, though it always has a selfadjoint
one in a larger complex Hilbert space (cf.\
\cite[Theorem 1 in Appendix I.2]{a-g}). This means
that each symmetric operator is subnormal. We say that
$A$ is {\em nonnegative} if $\is{Ah}{h} \Ge 0$ for all
$h \in \dz{A}$. Given two nonnegative selfadjoint
operators $C$ and $D$ in $\hh$, we write $C \preceq D$
if $\dz {D^{1/2}} \subseteq \dz{C^{1/2}}$ and
$\|C^{1/2} h\| \Le \|D^{1/2} h\|$ for all $h \in \dz
{D^{1/2}}$; note that $C \preceq D$ if and only if
$(D+x I)^{-1} \Le (C+x I)^{-1}$ for all real $x > 0$
or equivalently for some real $x > 0$ (cf.\
\cite[Theorem VI.2.21]{kato}). If $A$ is densely
defined and nonnegative, then there exist nonnegative
selfadjoint operators $B_{\mathrm K}$ and $B_{\mathrm
F}$ in $\hh$ that extends $A$ and such that
$B_{\mathrm K} \preceq B \preceq B_{\mathrm F}$ for
every nonnegative selfadjoint extension $B$ of $A$ in
$\hh$. The operators $B_{\mathrm K}$ and $B_{\mathrm
F}$ are called the {\em Krein} and the {\em
Friedrichs} extensions of $A$. We refer the reader to
\cite{b-s,weid} and
\cite{Cod-Sn,Se-St,Pr-Seb1,Pr-Seb2} for more
information on these subjects.

An operator $A$ in $\hh$ is called {\em paranormal} if
$\|Af\|^2 \Le \|f\|\|A^2 f\|$ for all $f \in
\dz{A^2}$. We say that an operator $A$ in $\hh$ is
{\em hyponormal} if $A$ is densely defined, $\dz{A}
\subseteq \dz{A^*}$ and $\|A^* f\| \Le \|Af\|$ for all
$f\in \dz{A}$. A densely defined operator $N$ in $\hh$
is said to be {\em normal} if $N$ is closed and $N^*N
= NN^*$ (or equivalently if and only if $N$ is closed
and both operators $N$ and $N^*$ are hyponormal, cf.\
\cite[Section 5.6]{weid}). A densely defined operator
$S$ in $\hh$ is called {\em subnormal} if there exists
a complex Hilbert space $\kk$ and a normal operator
$N$ in $\kk$ such that $\hh \subseteq \kk$ (isometric
embedding) and $Sh = Nh$ for all $h \in \dz S$. It is
well-known that normality implies subnormality,
subnormality implies hyponormality and hyponormality
implies paranormality, but none of these implications
can be reversed in general, i.e.,
   \begin{align*}
\{\textrm{normals}\} \varsubsetneq
\{\textrm{subnormals}\} \varsubsetneq
\{\textrm{hyponormals}\} \varsubsetneq
\{\textrm{paranormals}\}.
   \end{align*}
For details on this we refer the reader to
\cite{Con,Hal,furr,j-j-s} (see also
\cite{weid,b-s,StSz4,jj3,ot-sch,sz1} for the unbounded
case).
   \section{THE CLASSICAL MOMENT PROBLEM REVISITED}
   \subsection{\label{sec3}Indeterminate
moment problems} A sequence $\gams$ of real numbers is
said to be a {\em Stieltjes moment sequence} if there
exists a positive Borel measure $\mu$ on $\rbb_+$ such
that (from now on, we abbreviate $\int_{\rbb_+}$ to
$\int_0^\infty$)
   \begin{align*}
   \gamma_{n}=\int_0^\infty x^n \D\mu(x),\quad n\in
\zbb_+.
   \end{align*}
Call such $\mu$ an {\em \Sti-representing measure} of
the Stieltjes moment sequence $\gams$. A Stieltjes
moment sequence is said to be {\em \Sti-determinate}
if it has only one \Sti-representing measure;
otherwise, we call it {\em \Sti-indeterminate}. By the
Stieltjes theorem (cf.\ \cite[Theorem 6.2.5]{ber}), a
sequence $\gams \subseteq \rbb$ is a Stieltjes moment
sequence if and only if the sequences $\gams$ and
$\gamsj$ are positive definite (recall that a sequence
$\gams \subseteq \rbb$ is said to be {\em positive
definite} if $\sum_{k,l=0}^n \gamma_{k+l} \alpha_k
\overline{\alpha_l} \Ge 0$ for all $\alpha_0,\ldots,
\alpha_n \in \cbb$ and $n \in \zbb_+$). It is clear
that if $\gams$ is a Stieltjes moment sequence, then
so is $\gamsj$. The converse is easily seen to be
false (consider, e.g., the sequence
$\gams:=\{\gamma_0,1, 0, 0, \ldots\}$). Moreover, if a
Stieltjes moment sequence $\gams$ is
\Sti-indeterminate, then so is $\gamsj$ (see
\cite[Proposition 5.12]{sim}; see also Lemma
\ref{bext} below). The converse implication fails to
hold (cf.\ \cite[Corollary 4.21]{sim}; see also the
discussion below).

The following result has been established in
\cite{b-j-j-s} (see also \cite{wri} and \cite{sz} for
the question of backward extendibility of Hamburger
moment sequences).
   \begin{lem}[\mbox{\cite[Lemma 2.4.1]{b-j-j-s}}]
   \label{bext} Let $\gams$ be a Stieltjes moment
sequence and let $\gamma_{-1}$ be a positive real
number. Then the following are
equivalent\/\footnote{\;\label{foot}We adhere to the
convention that $\frac 1 0 := \infty$. Hence,
$\int_0^\infty \frac 1 x\D \mu(x) < \infty$ implies
$\mu(\{0\})=0$.}\/{\em :}
   \begin{enumerate}
   \item[(i)] $\gamsmj$ is a
Stieltjes moment sequence,
   \item[(ii)] there exists an \Sti-representing measure $\mu$
of $\gams$ such that
   \begin{align}   \label{bek}
\int_0^\infty \frac 1 x\D \mu(x) \Le \gamma_{-1}.
   \end{align}
   \end{enumerate}
Moreover, if {\em (i)} holds, then the mapping
$\mm_0(\gamma_{-1}) \ni \mu \to \nu_{\mu} \in
\mm_{-1}(\gamma_{-1})$ defined~ by
   \begin{align}   \label{nu}
\nu_{\mu}(\sigma) = \int_\sigma \frac 1 x\D \mu(x) +
\Big(\gamma_{-1} - \int_0^\infty \frac 1 x\D
\mu(x)\Big) \delta_0(\sigma), \quad \sigma \in
\borel{\rbb_+},
   \end{align}
is a bijection with the inverse $\mm_{-1}(\gamma_{-1})
\ni \nu \to \mu_{\nu} \in\mm_0(\gamma_{-1})$ given by
   \begin{align*}
\mu_{\nu} ( \sigma) = \int_\sigma x\D \nu (x), \quad
\sigma \in \borel{\rbb_+},
   \end{align*}
where $\mm_0(\gamma_{-1})$ is the set of all
\Sti-representing measures $\mu$ of $\gams$ such that
$\int_0^\infty \frac 1 x\D \mu(x) \Le \gamma_{-1}$,
and $\mm_{-1}(\gamma_{-1})$ is the set of all
\Sti-representing measures $\nu$ of $\gamsmj$. In
particular, $\nu_{\mu}(\{0\})=0$ if and only if
$\int_0^\infty \frac 1 x\D \mu(x)=\gamma_{-1}$.

If {\em (i)} holds and the sequence $\gams$ is
\Sti-determinate, then $\gamsmj$ is \Sti-determinate,
the unique \Sti-representing measure $\mu$ of $\gams$
satisfies the inequality $\int_0^\infty \frac 1 x \D
\mu(x) \Le \gamma_{-1}$, and $\nu_{\mu}$ is the unique
\Sti-representing measure of $\gamsmj$.
   \end{lem}
A sequence $\gams \subseteq \rbb$ is said to be a {\em
Hamburger moment sequence} if there exists a positive
Borel measure $\mu$ on $\rbb$ such that
   \begin{align*}
\gamma_{n}=\int_{-\infty}^\infty x^n \D\mu(x),\quad
n\in \zbb_+.
   \end{align*}
Call such $\mu$ an {\em \Ham-representing measure} of
the Hamburger moment sequence $\gams$. A Hamburger
moment sequence is said to be {\em \Ham-determinate}
if it has only one \Ham-representing measure;
otherwise, we call it {\em \Ham-indeterminate}. By the
Hamburger theorem (cf.\ \cite[Theorem 6.2.2]{ber}), a
sequence $\gams \subseteq \rbb$ is a Hamburger moment
sequence if and only if it is positive definite. It is
clear that if a Stieltjes moment sequence is
\Sti-indeterminate, then it is \Ham-indeterminate. The
reverse implication is not true in general (cf.\
\cite[page 96]{sim}).

Let $\gams$ be an \Ham-indeterminate Hamburger moment
sequence. By an {\em N-extremal} measure of $\gams$ we
mean an \Ham-representing measure $\mu$ of $\gams$ for
which the complex polynomials in one variable are
dense in $L^2(\mu)$. It is well-known that there is a
bijection $t \mapsto \mu_t$ between the set $\rbb \cup
\{\infty\}$ and the set of all N-extremal measures of
$\gams$ such that (cf.\ \cite[Remark, page 96]{sim})
   \begin{align} \label{Zenspi}
\int_0^\infty \frac{\D \mu_t(x)}{x} = t, \quad t \in
\rbb \cup \{\infty\}.
   \end{align}
The parametrization $t \mapsto \mu_t$ can be done as
follows (cf.\ \cite{sim}). Denote by $\pcal$ the ring
of all polynomials in one formal variable $X$ with
complex coefficients. Since $\gams$ is indeterminate,
there exists a unique inner product
$\is{\cdot}{\mbox{-}}$ on $\pcal$ such that
   \begin{align}    \label{xmxn}
\is{X^m}{X^n} = \gamma_{m+n}, \quad m,n\in \zbb_+.
   \end{align}
Let $\hh$ be the complex Hilbert space completion of
$(\pcal,\is{\cdot}{\mbox{-}})$. Since
$\is{Xp}{q}=\is{p}{Xq}$ for all $p, q\in \pcal$, we
deduce that there exists a unique symmetric operator
$A$ in $\hh$ such that $\dz{A}=\pcal$ and $A(p) = X
\cdot p$ for all $p\in \pcal$. Then clearly $\dz{A}$
is equal to the linear span of $\{A^ne\colon n \in
\zbb_+\}$ and, by \eqref{xmxn},
   \begin{align} \label{repr}
\gamma_n = \is{A^n e}{e}, \quad n \in \zbb_+ \quad
(e:=X^0).
   \end{align}
Hence, if $B$ is a von Neumann extension of $A$, then
$\mu_B(\cdot):= \is{E_B(\cdot)e}{e}$ is an
\Ham-representing measure of $\gams$, where $E_B$ is
the spectral measure of $B$. By the \Ham-indeterminacy
of $\gams$, the symmetric operator $A$ is not
essentially selfadjoint and its deficiency indices are
both equal to $1$, and thus there exists a bijection
$t \mapsto B_t$ between the set $\rbb \cup \{\infty\}$
and the set of all von Neumann extensions of $A$ such
that for every $t \in \rbb$, the spectrum of $B_t$
does not contain $0$ and $t=\is{B_t^{-1}e}{e}$, and
$0$ is an eigenvalue of $B_{\infty}$ (see
\cite[formulas (4.20)]{sim}
and\footnote{\;Unfortunately, the independence
assertion of \cite[Theorem 2.6]{sim}, saying that the
family $\{B_t\}_{t \in \rbb \cup \{\infty\}}$ is
independent of the choice of $\psi$, is not true (see
Appendix). Fortunately, the choice of $\psi$ made in
\cite[(4.20)]{sim} suits both the Hamburger and
Stieltjes moment problems.} \cite[Theorem 2.6]{sim}).
This immediately implies \eqref{Zenspi} with
$\mu_t(\cdot):= \is{E_{B_t}(\cdot)e}{e}$ for $t \in
\rbb\cup \{\infty\}$. It turns out that for every $t
\in \rbb \cup \{\infty\}$, $\mu_t$ is an N-extremal
measure of $\gams$ (and that there are no other
N-extremal measures), the closed support of $\mu_t$
(which coincides with the spectrum of $B_t$) has no
accumulation point in $\rbb$, and consequently it is
infinite and countable. Moreover, $\supp{\mu_s} \cap
\supp{\mu_t} = \varnothing$ for all $s,t \in \rbb \cup
\{\infty\}$ such that $s\neq t$, and
$\rbb=\bigcup_{t\in \rbb \cup \{\infty\}}
\supp{\mu_t}$, which means that the family
$\{\supp{\mu_t}\}_{t\in \rbb \cup \{\infty\}}$ forms a
partition of $\rbb$.

Now suppose that $\gams$ is an \Sti-indeterminate
Stieltjes moment sequence. Then $\gams$ is
\Ham-indeterminate. Let $(\hh,e,A)$ be as above. Then
$A$ is nonnegative (in fact $A - \alpha I$ is
nonnegative for some real $\alpha>0$) and it has many
nonnegative selfadjoint extensions in $\hh$. As a
consequence, the Krein extension $B_{\mathrm K}$ of
$A$ is different from the Friedrichs extension
$B_{\mathrm F}$ of $A$. It follows from \cite[Theorem
4.18]{sim} that $B_{\mathrm K} = B_{\infty}$ and
$B_{\mathrm F}=B_{t_0}$, where $t_0 = \is{B_{\mathrm
F}^{-1}e}{e} \in (0, \infty)$, and\footnote{\;See also
\cite[Theorem 5.2]{ber1} for a Nevanlinna type
parametrization of solutions of an \Sti-indeterminate
Stieltjes moment sequence. Both parameterizations are
equivalent.}
   \begin{align} \label{exp}
\forall t \in \rbb \cup \{\infty\}\colon \;
\supp{\mu_t} \subseteq [0,\infty) \iff t \in
[t_0,\infty) \cup \{\infty\}.
   \end{align}
In other words, $\{\mu_t\}_{t\in [t_0,\infty) \cup
\{\infty\}}$ are the only N-extremal measures of
$\gams$ which are simultaneously \Sti-representing
measures of $\gams$. Call the N-extremal measures
$\muk(\cdot) := \is{E_{B_\infty}(\cdot)e}{e}$ and
$\muf(\cdot) := \is{E_{B_{t_0}}(\cdot)e}{e}$ the {\em
Krein} and the {\em Friedrichs} measures of $\gams$,
respectively. Note that $\muk=\mu_{\infty}$ and
$\muf=\mu_{t_0}$. Arguing as in the proof of
\cite[Proposition 3.1]{sim}, we deduce that
$\min(\supp{\mu_t}) < \min(\supp{\muf})$ for all $t
\in (t_0,\infty)\cup\{\infty\}$. Hence, by the
preceding paragraph and \eqref{exp}, we have
   \begin{align*}
\text{$0 \in \supp{\muk}$ and $0 < \min(\supp{\mu_t})
< \min(\supp{\muf})$ for all $t \in (t_0,\infty)$.}
   \end{align*}
This in turn implies that
   \begin{align} \label{exp1}
\text{$0 < \int_0^\infty \frac{1}{x^n} \D\mu_t(x) <
\infty$ for all $n \in \nbb$ and $t \in
[t_0,\infty)$.}
   \end{align}
In particular $0 < \int_0^\infty \frac{1}{x^n}
\D\muf(x) < \infty$ for all $n \in \nbb$.
   \subsection{\label{sec4}Krein and Friedrichs measures}
Now we state some crucial inequalities for the Krein
and Friedrichs measures.
   \begin{thm}[\mbox{\cite[Theorem 4.19 and
Corollary 4.20]{sim}}] \label{1} Let $\gams$ be an
\Sti-indeterminate Stieltjes moment sequence and let
$\muk$, $\muf$ be the corresponding Krein and
Friedrichs measures. If $\rho$ is an \Sti-representing
measure of $\gams$ such that $\rho\neq \muf$, then
   \begin{align}  \label{sim1}
\int_0^\infty \frac{\D\muf(x)}{x+y} < \int_0^\infty
\frac{\D\rho(x)}{x+y} \Le \int_0^\infty
\frac{\D\muk(x)}{x+y}, \quad y \in [0,\infty).
   \end{align}
   \end{thm}
   \begin{cor}
Let $\boldsymbol{\gamma}=\gams$ be an
\Sti-indeterminate Stieltjes moment sequence and let
$\EuScript M^{\mathrm S}(\boldsymbol{\gamma})$ be the
set of all its \Sti-representing measures. Then the
Friedrichs measure $\muf$ of $\boldsymbol{\gamma}$ is
a unique measure $\rho\in \EuScript M^{\mathrm
S}(\boldsymbol{\gamma})$ such that
   \begin{align*}
\int_0^\infty \frac{1}{x} \D\rho(x) = \min
\Big\{\int_0^\infty \frac{1}{x} \D\sigma (x)\colon
\sigma \in \EuScript M^{\mathrm
S}(\boldsymbol{\gamma}) \Big\}.
   \end{align*}
   \end{cor}
We will show that the right-hand inequality in
\eqref{sim1} is in fact strict for all real $y>0$ (but
not for $y=0$ as explained just after the proof of
Proposition~ \ref{1+}). This is an answer to a
question raised by C. Berg \cite{Berg-pc}.
   \begin{pro} \label{1+} Let $\gams$ be an
\Sti-indeterminate Stieltjes moment sequence and let
$\muk$ be its Krein measure. If $\rho$ is an
\Sti-representing measure of $\gams$ such that
$\rho\neq \muk$, then
   \begin{align}   \label{sim2}
\int_0^\infty \frac{\D\rho(x)}{x+y} < \int_0^\infty
\frac{\D\muk(x)}{x+y}, \quad y \in (0,\infty).
   \end{align}
   \end{pro}
   \begin{proof}
It follows from \cite[Theorem 4.18]{sim} that there
are entire functions $A,B,C,D$ (determined by the
sequence $\gams$) such that for all $t \in [t_0,
\infty) \cup \{\infty\}$,
   \begin{align} \label{Nevan}
F(-y)(t):= - \frac{C(-y)t + A(-y)}{D(-y)t + B(-y)} =
\int_0^\infty \frac{\D\mu_t(x)}{x+y}, \quad y \in
(0,\infty),
   \end{align}
where the middle term in \eqref{Nevan} is understood
as $- \frac{C(-y)}{D(-y)}$ for $t=\infty$. Since $A,
B, C, D$ take real values on the real line and $AD -
BC\equiv 1$ (cf.\ \cite[Theorem 4.8(iii)]{sim}), we
deduce that the derivative of $F(-y)(\cdot)$ is
positive on $[t_0,\infty)$, and thus the map
$F(-y)(\cdot)$ is strictly increasing on
$[t_0,\infty)$. Then for all $t \in [t_0,\infty)$,
   \begin{multline} \label{nier}
F(-y)(t) = - \frac{C(-y)t + A(-y)}{D(-y)t + B(-y)}
\underset{(t \to \infty)}{\nearrow} -
\frac{C(-y)}{D(-y)} = F(-y)(\infty).
   \end{multline}

If the measure $\rho$ is N-extremal, then by our
assumption and \cite[Theorem 4.18]{sim} there is $t
\in [t_0,\infty)$ such that $\rho= \mu_t$. Then, by
\eqref{nier}, we have
   \begin{align*}
\int_0^\infty \frac{\D\mu_t(x)}{x+y}
\overset{\eqref{Nevan}}= F(-y)(t) < F(-y)(\infty)
\overset{\eqref{Nevan}}= \int_0^\infty
\frac{\D\muk(x)}{x+y}, \quad y \in (0,\infty).
   \end{align*}
If $\rho$ is not N-extremal, then, again by
\cite[Theorem 4.18]{sim}, there is a non-constant Pick
function $\varPhi\colon \cbb \setminus [0,\infty) \to
\cbb$ such that $\varPhi(-y) \in [t_0,\infty)$ for all
$y \in (0,\infty)$, and
   \begin{align} \label{podstaw}
\int_0^\infty \frac{\D\rho(x)}{x-z} = -
\frac{C(z)\varPhi(z) + A(z)}{D(z)\varPhi(z) + B(z)},
\quad z \in \cbb \setminus [0,\infty).
   \end{align}
Hence, substituting $z=-y$ into \eqref{podstaw}, we
get
   \begin{align*}
\int_0^\infty \frac{\D\rho(x)}{x+y} =
F(-y)(\varPhi(-y)) \overset{\eqref{nier}} <
F(-y)(\infty) \overset{\eqref{Nevan}}= \int_0^\infty
\frac{\D\muk(x)}{x+y}, \quad y \in (0,\infty).
   \end{align*}
This completes the proof.
   \end{proof}
   \begin{rem}
Let $\gams$ be any \Sti-indeterminate Stieltjes moment
sequence. Fix $\alpha \in (0,1)$ and set
$\rho_{\alpha}= \alpha \muk + (1-\alpha) \muf$, where
$\muk$ and $\muf$ are the Krein and the Friedrichs
measures of $\gams$. Then $\rho_{\alpha}$ is an
\Sti-representing measure of $\gams$ such that
$\rho_{\alpha} \neq \muk$, $\rho_{\alpha}$ is not
N-extremal and, because $0$ is an atom of $\muk$,
   \begin{align*}
\int_0^\infty \frac{\D\rho_{\alpha}(x)}{x} =
\int_0^\infty \frac{\D\muk(x)}{x} = \infty.
   \end{align*}
In other words, the strict inequality in \eqref{sim2}
may turn into equality when $y=0$. This is never the
case for N-extremal measure $\rho$ (apply
\eqref{Zenspi}).
   \end{rem}
Before stating the next result, we prove a lemma which
is of some independent interest.
   \begin{lem} \label{oin}
If $\gams$ is an \Sti-determinate Stieltjes moment
sequence whose \Sti-representing measure $\tau$ has
the property that $\tau(\{0\})=0$, then $\gams$ is
\Ham-determinate.
   \end{lem}
   \begin{proof}
Suppose that, contrary to our claim, $\gams$ is
\Ham-indeterminate. Then the operator $A$ attached to
$\gams$ via \eqref{repr} is not essentially
selfadjoint. Since $\gams$ is \Sti-determinate, we
deduce from \cite[Theorem 2]{sim} (see also
\cite[Theorem 5]{c-s-sz}) that $A$ has a unique
nonnegative selfadjoint extension in $\hh$ which is
evidently the Friedrichs extension $B_{\mathrm F}$ of
$A$. Hence, by \cite[Proposition 3.1]{sim}, $0$ is an
eigenvalue of $B_{\mathrm F}$. Denote by $E$ the
spectral measure of $B_{\mathrm F}$. Then clearly
$\mu(\cdot):=\is{E(\cdot)e}{e}$ is an N-extremal
measure of $\gams$. Since the closed support of any
N-extremal measure has no accumulation point in $\rbb$
and $\supp\mu$ coincides with the spectrum of
$B_{\mathrm F}$ (see \cite[Theorem 5]{c-s-sz} and also
\cite[Theorem 5]{StSz4}), we deduce that $\mu$ is an
\Sti-representing measure of $\gams$ and $0$ is an
atom of $\mu$. By the \Sti-determinacy of $\gams$, we
have $\tau=\mu$, which implies that
$\tau(\{0\})\neq0$, a contradiction. This completes
the proof.
   \end{proof}
Note that if $\gams$ is an \Ham-determinate Stieltjes
moment sequence, then its unique \Ham-representing
measure may have an atom at $0$ (any compactly
supported finite positive Borel measure on
$[0,\infty)$ with an atom at $0$ is an
\Ham-representing measure of such a sequence). This
means that the converse of the implication in Lemma
\ref{oin} does not hold in general.

The following characterization of the \Ham-determinacy
of a borderline backward extension of an
\Sti-indeterminate Stieltjes moment sequence will be
used in the proof of Theorem \ref{wkwsti}. Let us
mention that the implication (ii)$\Rightarrow$(i) and
the ``moreover'' part of Theorem \ref{2} below has
appeared in \cite[Corollary 4.21]{sim}. We include
their proofs to keep the exposition self-contained.
   \begin{thm} \label{2}
Let $\gams$ be an \Sti-indeterminate Stieltjes moment
sequence, $\muf$ be its Friedrichs measure and
$\gamma_{-1}$ be a nonnegative real number. Then the
following two conditions are equivalent{\em :}
   \begin{enumerate}
   \item[(i)] $\gamsmj$ is an \Sti-determinate Stieltjes moment
sequence,
   \item[(ii)] $\gamma_{-1}=\int_0^\infty
\frac{\D\muf(x)}{x}$.
   \end{enumerate}
Moreover, if any of the above equivalent conditions
holds, then $\gamsmj$ is \Ham-determinate.
   \end{thm}
   \begin{proof}
Let $\{\mu_t\}_{t\in\rbb \cup \{\infty\}}$ be the
parametrization of N-extremal measures of $\gams$
given by \eqref{Zenspi}.

   (i)$\Rightarrow$(ii) Note that $\gamma_{-1}
> 0$ (otherwise  $\gamma_n=0$ for
all $n\in \zbb_+$). By the \Sti-determinacy of
$\gamsmj$ and Lemma \ref{bext}, there is a unique
\Sti-representing measure $\rho$ of $\gams$ such that
$\int_0^\infty \frac 1 x \D \rho(x) \Le \gamma_{-1}$.
In view of Theorem \ref{1}, we~ have
   \begin{align} \label{ququ}
t_0\overset{\eqref{Zenspi}}=\int_0^\infty \frac 1 x \D
\muf(x) \Le \int_0^\infty \frac 1 x \D \rho(x) \Le
\gamma_{-1} \overset{\eqref{Zenspi}}= \int_0^\infty
\frac{1}{x} \D \mu_{\gamma_{-1}}(x),
   \end{align}
which, by \eqref{exp}, implies that
$\mu_{\gamma_{-1}}$ is an \Sti-representing measure of
$\gams$. Since, by \eqref{ququ}, $\muf$ and
$\mu_{\gamma_{-1}}$ satisfy inequality \eqref{bek}, we
conclude that $\muf=\rho=\mu_{\gamma_{-1}}$. This
gives (ii).

   (ii)$\Rightarrow$(i) If $\gamsmj$ were not
\Sti-determinate, then by Lemma \ref{bext}, there
would exist an \Sti-representing measure $\rho$ of
$\gams$ such that $\rho\neq \muf$ and
   \begin{align*}
\int_0^\infty \frac{1}{x} \D \rho (x) \Le
\int_0^\infty \frac{1}{x} \D \muf(x),
   \end{align*}
which would contradict \eqref{sim1}.

If $\gamsmj$ is \Sti-determinate, then by \eqref{nu}
with $\mu=\muf$ we see that $\D \tau(x):=\frac 1 x\D
\muf(x)$ is an \Sti-representing measure of $\gamsmj$
such that $0\notin \supp{\tau}$ (because $0\notin
\supp{\muf}$). Hence the ``moreover'' part follows
from Lemma \ref{oin}.
   \end{proof}
We are now ready to state a result which is the main
tool for constructing an operator with properties
mentioned in the title of the paper.
   \begin{thm} \label{wkwsti}
Suppose that $\gams$ is an \Sti-indeterminate
Stieltjes moment sequence, $\muf$ is its Friedrichs
measure and $\gamma_{-1}$ is a nonnegative real
number. Then the following assertions hold.
   \begin{enumerate}
   \item[(i)] If $\gamma_{-1} < \int_0^\infty \frac 1 x \D
\muf(x)$, then $\gamsmj$ is not a Stieltjes moment
sequence.
   \item[(ii)] If $\gamma_{-1} =  \int_0^\infty \frac 1 x \D
\muf(x)$, then $\gamsmj$ is an \Ham-determinate
Stieltjes moment sequence.
   \item[(iii)] If $\gamma_{-1} >  \int_0^\infty \frac 1 x \D
\muf(x)$, then $\gamsmj$ is an \Sti-indeterminate
Stieltjes moment sequence.
   \end{enumerate}
   \end{thm}
   \begin{proof}
Assertions (i) and (ii) follow from Lemma \ref{bext}
and Theorems \ref{1} and \ref{2}.

(iii) By Lemma \ref{bext}, the sequence $\gamsmj$ is a
Stieltjes moment sequence. In view of \eqref{Zenspi}
and \eqref{exp}, the measures $\muf$ and
$\mu_{\gamma_{-1}}$ are two distinct \Sti-representing
measures of $\gams$ which satisfy \eqref{bek}. Hence,
by Lemma \ref{bext}, $\gamsmj$ is an
\Sti-indeterminate Stieltjes moment sequence.
   \end{proof}
   \begin{cor}
Let $\gams$ be an \Sti-indeterminate Stieltjes moment
sequence and let $\muf$ be the Friedrichs measure of
$\gams$. Then
   \begin{align} \label{A2}
\int_0^\infty \frac{1}{x} \D\muf(x) =
\min\Big\{\gamma_{-1} \in (0,\infty)\colon \forall
n\Ge 0 \;\; \det[\gamma_{i+j-1}]_{i,j=0}^n > 0\Big\}.
   \end{align}
   \end{cor}
   \begin{proof}
Set $t_0=\int_0^\infty \frac{1}{x} \D\muf(x)$. It
follows from Theorem \ref{wkwsti} that
   \begin{align*}
t_0 = \min\Big\{\gamma_{-1} \in (0,\infty) \colon
\gamsmj \text{ is a Stieltjes moment sequence}\Big\}.
   \end{align*}
Applying the Stieltjes and Hamburger theorems (cf.\
\cite[Theorems 6.2.5 and 6.2.2]{ber}) and using the
fact that $\gams$ is a Hamburger moment sequence, we
deduce that
   \begin{align*}
t_0 = \min\Big\{\gamma_{-1} \in (0,\infty) \colon
\gamsmj \text{ is a Hamburger moment sequence}\Big\}.
   \end{align*}
This equality, when combined with \cite[Theorem
1.2]{sh-tam} and the fact that $\gamsmj$ can never
have a finitely supported \Ham-representing measure
complete the proof.
   \end{proof}
   \subsection{\label{sec5}Peculiar Stieltjes moment sequences}
Our main objective here is to construct
\Sti-indeterminate Stieltjes moment sequences with
specific properties that will be used later to build
non-hyponormal operators generating Stieltjes moment
sequences.
   \begin{exa} \label{e1}
Fix $\kappa \in \zbb_+ \sqcup \{\infty\}$. We will
indicate a system $\{\gamma_n\}_{n=-\kappa}^\infty$ of
positive real numbers which has the following
properties:
   \begin{enumerate}
   \item[(i)] $\gamma_0=1$,
   \item[(ii)] there exists a positive
Borel measure $\nu$ on $(0,\infty)$ such that
   \begin{align*}
\gamma_n = \int_0^\infty x^n \D \nu(x), \quad n \in
\zbb, \, n \Ge -\kappa,
   \end{align*}
   \item[(iii)] $\gamsj$ is an \Sti-indeterminate
Stieltjes moment sequence,
   \item[(iv)] there exists an \Sti-representing
measure $\rho$ of $\gamsj$ such that
   \begin{gather}  \label{ineq00}
\text{$\supp{\rho}$ has no accumulation point in
$(0,\infty)$,}
   \\  \label{ineq0}
0<\int_0^\infty \frac{1}{x^n} \D \rho(x) < \infty,
\quad n =1, \ldots, \kappa+1,
   \end{gather}
and
   \begin{gather} \label{ineq1}
\int_0^\infty \frac{1}{x} \D \rho(x)
> 1.
   \end{gather}
   \end{enumerate}
What is more, we can always construct a system
$\{\gamma_n\}_{n=-\kappa}^\infty$ of positive real
numbers which satisfies the conditions (i) to (iv) and
which has the property that the sequence $\gams$ is
either \Ham-determinate or \Sti-indeterminate
according to our needs.

For this purpose, we fix $q\in (0,1)$ and define
   \begin{align*}
\zeta_n=q^{-\frac {1}{2} n^2}, \quad n\in \zbb.
   \end{align*}
It is easily seen that for every $\theta \in [-1,1]$,
   \begin{align*}
\zeta_n = \int_0^\infty x^{n} \omega_{\theta}(x)\D x,
\quad n \in \zbb,
   \end{align*}
where the density function $\omega_{\theta}$ is given
by
   \begin{align*}
\omega_{\theta}(x)=\frac{1}{\sqrt{2 \pi} \, \sigma} \,
x^{-1} \exp\Big(- \frac{(\log x)^2}{2
\sigma^2}\Big)\left(1 + \theta
\sin\Big(\frac{2\pi}{\sigma^2} \log x \Big)\right),
\quad x \in (0,\infty),
   \end{align*}
with $\sigma = \sqrt{-\log q}$. This means that for
every $l\in \zbb$, the sequence
$\{\zeta_{n+l}\}_{n=0}^\infty$ is an
\Sti-indeterminate Stieltjes moment sequence. This is
a famous example due to Stieltjes (cf.\ \cite{sti}).
It was noticed much later by Chihara \cite{chi} and
Leipnik \cite{lei} (see also \cite{ber2}) that for
every $a\in (0,\infty)$, the Borel probability measure
$\lambda_a$ defined by
   \begin{align} \label{nosnik}
\lambda_a = \frac{1}{L(a)} \sum_{k=-\infty}^{\infty}
a^k q^{\frac 12 {k^2}} \delta_{aq^k}, \quad L(a)=
\sum_{k=-\infty}^{\infty} a^k q^{\frac 12 k^2},
   \end{align}
solves the moment problem
   \begin{align}  \label{oldsti2}
\zeta_n = \int_0^\infty x^n \D \lambda_a(x), \quad n
\in \zbb.
   \end{align}
Therefore, for every fixed $l \in \zbb$, the
absolutely continuous measures $x^l \omega_{\theta}(x)
\D x$, $\theta \in [-1,1]$, and the pure point
measures $x^l \D \lambda_a(x)$, $a\in (0,\infty)$, are
\Sti-representing measures of
$\{\zeta_{n+l}\}_{n=0}^\infty$. Since $0$ is an
accumulation point of the closed support of each of
these measures, we conclude that neither of them is
N-extremal.

Let $\{\mu_t\}_{t \in \rbb \cup \{0\}}$ be the set of
all N-extremal measures of
$\{\zeta_{n}\}_{n=0}^\infty$ (cf.\ \eqref{Zenspi}) and
let $\muf$ be the Friedrichs measure of
$\{\zeta_{n}\}_{n=0}^\infty$. Set $t_0 = \int_0^\infty
\frac {1}{x} \D\muf(x)$. Take $t\in [t_0,\infty)$ and
define the system $\{\gamma_n(t)\}_{n=-\kappa}^\infty$
by
   \begin{align*}
\gamma_n(t) =
   \begin{cases}
t^{-1} \int_0^\infty x^{n-1} \D\mu_{t}(x) & \text{ if
} -\kappa \Le n \Le 0,
   \\[1ex]
t^{-1} \zeta_{n-1} & \text{ if } n \Ge 1.
   \end{cases}
   \end{align*}
By \eqref{exp1}, the above definition is correct. It
is clear that the system
$\{\gamma_n(t)\}_{n=-\kappa}^\infty$ satisfies the
conditions (i) and (ii) with a measure $\nu$ given by
$\D\nu(x)=t^{-1}\frac{1}{x}\D\mu_{t}(x)$. Since
$\gamma_{n+1}(t) = t^{-1} \zeta_n$ for all $n \in
\zbb_+$, we see that the system
$\{\gamma_n(t)\}_{n=-\kappa}^\infty$ satisfies the
condition (iii) and that for every $s \in (t,
\infty)$, $\rho_s := t^{-1}\mu_s$ is an
\Sti-representing measure of $\gamstj{t}$ which
satisfies \eqref{ineq00} and \eqref{ineq0} (see
\eqref{exp1}). Moreover, we have
   \begin{align*}
\int_0^\infty \frac 1 x \D \rho_s(x) = t^{-1}
\int_0^\infty \frac 1 x \D \mu_{s}(x)
\overset{\eqref{Zenspi}}= t^{-1}s > 1, \quad s\in (t,
\infty),
   \end{align*}
which means that $\rho_s$ satisfies \eqref{ineq1} for
every $s\in (t, \infty)$. It follows from
\eqref{Zenspi} and Theorem \ref{wkwsti} that the
Stieltjes moment sequence $\gamst{t}$ is
\Ham-determinate for $t=t_0$ and \Sti-indeterminate
for $t\in (t_0,\infty)$.

Since the closed supports of the measures $\rho_s$, $s
\in (t,\infty)$, are not explicitly known, we will
provide other examples of measures satisfying the
conditions \eqref{ineq00}, \eqref{ineq0} and
\eqref{ineq1}, the closed supports of which are
precisely given. According to Theorem \ref{1} and the
fact that $\lambda_a$ is not N-extremal, we have
   \begin{align*}
t_0 = \int_0^\infty \frac 1 x \D \muf(x) <
\int_0^\infty \frac 1 x \D \lambda_a (x)
\overset{\eqref{oldsti2}}= \zeta_{-1}, \quad a \in
(0,\infty),
   \end{align*}
which means that $[t_0,\zeta_{-1}) \neq \varnothing$.
Take $t \in [t_0,\zeta_{-1})$ and set $\tilde
\rho_a=\frac{1}{t}\lambda_a$ for $a \in (0,\infty)$.
Using \eqref{oldsti2}, we can easily verify that for
every $a \in (0,\infty)$, $\tilde \rho_a$ is an
\Sti-representing measure of $\gamstj{t}$ which
satisfies \eqref{ineq00}, \eqref{ineq0} and
\eqref{ineq1}. By \eqref{nosnik}, $\supp{\tilde
\rho_a} = \{a q^k\colon k\in \zbb\} \cup \{0\}$ for
every $a \in (0,\infty)$.
   \end{exa}
Note that the constant $t_0$ which plays an essential
role in Example \ref{e1} can be estimated by using
\eqref{A2}.
   \section{RELATING MOMENTS TO DIRECTED TREES}
   \subsection{\label{sec6}Weighted shifts on directed trees}
Let $\tcal=(V,E)$ be a directed tree ($V$ and $E$
stand for the sets of vertices and edges of $\tcal$,
respectively). If $\tcal$ has a root, which will
always be denoted by $\koo$, then we write
$V^\circ:=V\setminus \{\koo\}$; otherwise, we put
$V^\circ = V$. Set
   \begin{align*}
\dzi u = \{v\in V\colon (u,v)\in E\}, \quad u \in V.
   \end{align*}
A member of $\dzi u$ is called a {\em child} (or {\em
successor}) of $u$. For every vertex $u \in V^\circ$
there exists a unique vertex, denoted by $\pa u$, such
that $(\pa u,u)\in E$. The correspondence $u \mapsto
\pa u$ is a partial function from $V$ to $V$. For an
integer $n \Ge 1$, the $n$-fold composition of the
partial function $\paa$ with itself will be denoted by
$\paa^n$. Let $\paa^0$ stand for the identity map on
$V$. We call $\tcal$ {\em leafless} if $V = V^\prime$,
where $V^\prime:=\{u \in V \colon \dzi u \neq
\varnothing\}$. It is clear that every leafless
directed tree is infinite. A vertex $u \in V$ is said
to be a {\em branching vertex} of $\tcal$ if $\dzi{u}$
consists of at least two vertices. If $W \subseteq V$,
we put $\dzi W = \bigcup_{v \in W} \dzi v$ and $\des W
= \bigcup_{n=0}^\infty \dzin n W$, where $\dzin{0}{W}
= W$ and $\dzin{n+1}{W} = \dzi{\dzin{n}{W}}$ for all
integers $n\Ge 0$. For $u \in V$, we set $\dzin n
u=\dzin n {\{u\}}$ and $\des{u}=\des{\{u\}}$. It
follows from \cite[Proposition 2.1.2]{j-j-s} and
\cite[Proposition 2.2.1]{b-j-j-s} that
   \begin{align} \label{roz}
V^\circ & = \bigsqcup_{u\in V} \dzi u,
   \\
\label{dzinn2} \dzin{n+1}{u} & = \bigsqcup_{v \in
\dzi{u}} \dzin{n}{v}, \quad n \in \zbb_+,\, u \in V,
   \end{align}
where the symbol $\bigsqcup$ is reserved to denote
pairwise disjoint union of sets.

Given a directed tree $\tcal$, we tacitly assume that
$V$ and $E$ stand for the sets of vertices and edges
of $\tcal$, respectively. Denote by $\ell^2(V)$ the
complex Hilbert space of all square summable complex
functions on $V$ with the standard inner product. For
$u \in V$, we define $e_u$ to be the characteristic
function of the one-point set $\{u\}$. The family
$\{e_u\}_{u\in V}$ is an orthonormal basis of
$\ell^2(V)$. We write $\escr$ for the linear span of
the set $\{e_u\colon u \in V\}$.

Given $\lambdab = \{\lambda_v\}_{v \in V^\circ}
\subseteq \cbb$, we define the operator $\slam$ in
$\ell^2(V)$ by
   \begin{align*}
   \begin{aligned}
\dz {\slam} & = \{f \in \ell^2(V) \colon
\varLambda_\tcal f \in \ell^2(V)\},
   \\
\slam f & = \varLambda_\tcal f, \quad f \in \dz
{\slam},
   \end{aligned}
   \end{align*}
where $\varLambda_\tcal$ is the map defined on
functions $f\colon V \to \cbb$ via
   \begin{align} \label{lamtauf}
(\varLambda_\tcal f) (v) =
   \begin{cases}
\lambda_v \cdot f\big(\pa v\big) & \text{ if } v\in
V^\circ,
   \\
0 & \text{ if } v=\koo.
   \end{cases}
   \end{align}
The operator $\slam$ is called a {\em weighted shift}
on the directed tree $\tcal$ with weights
$\lambdab=\{\lambda_v\}_{v \in V^\circ}$. Combining
Propositions 3.1.2, 3.1.3\,(iii) and 3.1.7 of
\cite{j-j-s}, we get the ensuing properties of $\slam$
(from now on, we adopt the convention that
$\sum_{v\in\varnothing} x_v=0$).
   \begin{pro}\label{bas}
Let $\slam$ be a weighted shift on a directed tree
$\tcal$ with weights $\lambdab = \{\lambda_v\}_{v \in
V^\circ}$. Then the following assertions hold.
   \begin{enumerate}
   \item[(i)] $\slam$ is closed.
   \item[(ii)] $e_u$ is in $\dz{\slam}$ if and only if
$\sum_{v\in\dzi u} |\lambda_v|^2 < \infty$; if $e_u
\in \dz{\slam}$, then
   \begin{align*}
\slam e_u = \sum_{v\in\dzi u} \lambda_v e_v \quad
\text{and} \quad \|\slam e_u\|^2 = \sum_{v\in\dzi u}
|\lambda_v|^2.
   \end{align*}
   \item[(iii)] $\slam$ is injective if and only if
$\tcal$ is leafless and $\sum_{v\in\dzi u}
|\lambda_v|^2 > 0$ for every $u\in V$.
   \end{enumerate}
   \end{pro}
Let us now recall a characterization of hyponormality
of weighted shifts on leafless directed trees with
nonzero weights.
    \begin{thm}[\mbox{\cite[Theorem
5.1.2 and Remark 5.1.5]{j-j-s}}] \label{hyp} Let
$\slam$ be a densely defined weighted shift on a
leafless directed tree $\tcal$ with nonzero weights
$\lambdab = \{\lambda_v\}_{v \in V^\circ}$. Then
$\slam$ is hyponormal if and only if
    \begin{gather} \label{hypon}
\sum_{v \in \dzi{u}} \frac {|\lambda_v|^2}{\|\slam
e_v\|^2} \Le 1, \quad u \in V.
    \end{gather}
    \end{thm}
The following lemma relates representing measures of
Stieltjes moment sequences induced by basic vectors
coming from the parent and its children. Inequality
\eqref{con2} below will be referred to as the {\em
consistency condition} at $u$.
   \begin{lem}[\mbox{\cite[Lemma 4.2.3]{b-j-j-s}}]
Let $\slam$ be a weighted shift on a directed tree
$\tcal$ with weights $\lambdab=\{\lambda_v\}_{v \in
V^\circ}$ such that $\escr \subseteq \dzn{\slam}$. Let
$u \in V^\prime$. Suppose that for every $v \in \dzi
u$ the sequence $\{\|\slam^n e_v\|^2\}_{n=0}^\infty$
is a Stieltjes moment sequence with a representing
measure $\mu_v$. Consider the following two
conditions\footnote{\;\;We adhere to the standard
convention that $0 \cdot \infty = 0$; see also
footnote \ref{foot}.}{\em :}
   \begin{gather} \label{con1}
\text{$\{\|\slam^n e_u\|^2\}_{n=0}^\infty$ is a
Stieltjes moment sequence,}
   \\  \label{con2}
\sum_{v \in \dzi{u}} |\lambda_v|^2 \int_0^\infty \frac
1 x\, \D \mu_v(x) \Le 1.
   \end{gather}
   Then the following assertions are valid.
   \begin{enumerate}
   \item[(i)] If  \eqref{con2}  holds, then so does
\eqref{con1} and the positive Borel measure $\mu_u$ on
$\rbb_+$ defined by
   \begin{align} \label{muu+}
\mu_u(\sigma) = \sum_{v \in \dzi u} |\lambda_v|^2
\int_\sigma \frac 1 x \D \mu_v(x) + \varepsilon_u
\delta_0(\sigma), \quad \sigma \in \borel{\rbb_+},
      \end{align}
with
   \begin{align} \label{muu++}
\varepsilon_u=1 - \sum_{v \in \dzi u} |\lambda_v|^2
\int_0^\infty \frac 1 x \D \mu_v(x),
   \end{align}
is a representing measure of $\{\|\slam^n
e_u\|^2\}_{n=0}^\infty$.
   \item[(ii)] If \eqref{con1} holds and $\{\|\slam^{n+1}
e_u\|^2\}_{n=0}^\infty$ is determinate, then
\eqref{con2} holds, the Stieltjes moment sequence
$\{\|\slam^n e_u\|^2\}_{n=0}^\infty$ is determinate
and its unique representing measure $\mu_u$ is given
by \eqref{muu+} and \eqref{muu++}.
   \end{enumerate}
   \end{lem}
   \subsection{\label{sec7}Generating Stieltjes
moments on directed trees} We begin by recalling the
action of powers of $\slam$ on basic vectors $e_u$, $u
\in V$.
   \begin{lem}[\mbox{\cite[Lemma 4.1.1]{b-j-j-s}}]
\label{lem4} Let $\slam $ be a weighted shift on a
directed tree $\tcal$ with weights $\lambdab =
\{\lambda_v\}_{v \in V^\circ}$. Then the following
assertions hold for all $u \in V$ and $n \in \zbb_+$.
   \begin{enumerate}
   \item[(i)] $e_u \in \dz{\slam^n}$ if and only if
$\sum_{v \in \dzin{m}{u}} |\lambda_{u\mid v}|^2 <
\infty$ for all integers $m$ such that $1 \Le m \Le
n$.
   \item[(ii)] If $e_u \in \dz{\slam^n}$, then
   \begin{align}  \label{pow}
\slam^n e_u & = \sum_{v \in \dzin{n}{u}}
\lambda_{u\mid v} \, e_v,
   \\     \label{brak}
\|\slam^n e_u\|^2 & = \sum_{v \in \dzin{n}{u}}
|\lambda_{u\mid v}|^2,
   \end{align}
   \end{enumerate}
where
   \begin{align}    \label{luv}
\lambda_{u\mid v} =
   \begin{cases}
1 & \text{ if } v=u,
   \\
\prod_{j=0}^{n-1} \lambda_{\paa^{j}(v)} & \text{ if }
v \in \dzin{n}{u}, \, n \Ge 1.
   \end{cases}
   \end{align}
   \end{lem}
One can deduce from \eqref{luv} that
   \begin{align} \label{recfor2}
\lambda_{\pa v\mid w} & = \lambda_v \lambda_{v\mid w},
\quad v \in V^\circ, \, w\in \des v.
   \end{align}
The above lemma enables us to describe the powers of
$\slam$. Below we write $\sumo$ for the sum of a
series whose terms are mutually orthogonal.
   \begin{thm} \label{dn}
Let $\slam$ be a weighted shift on a directed tree
$\tcal$ with weights $\lambdab = \{\lambda_v\}_{v\in
V^\circ}$. Then the following assertions hold for any
$n \in \zbb_+$.
   \begin{enumerate}
   \item[(i)] A function $f \colon V \to \cbb$ belongs
to $\dz{\slam^n}$ if and only if
   \begin{align} \label{series2}
\sum_{u\in V} |f(u)|^2 \Bigg(\sum_{j=0}^n \, \sum_{v
\in \dzin{j}{u}} |\lambda_{u\mid v}|^2\Bigg) < \infty,
   \end{align}
with the usual convention that $0 \cdot \infty = 0$.
   \item[(ii)] If $f \in \dz{\slam^n}$, then
$e_u \in \dz{\slam^n}$ for every $u \in V$ such that
$f(u)\neq 0$, and
   \begin{align}  \label{aaa}
\slam^n f & = \underset{u\in V:\, f(u) \neq 0}{\sumo}
f(u) \slam^n e_u, \quad f \in \dz{\slam^n},
   \\ \label{bbb}
\|\slam^n f\|^2 &= \sum_{u\in V:\, f(u) \neq 0}
|f(u)|^2 \|\slam^n e_u\|^2, \quad f \in \dz{\slam^n}.
   \end{align}
   \item[(iii)] If $\escr\subseteq \dz{\slam^n}$, then
$\escr$ is a core of $\slam^n$.
   \item[(iv)] $\slam^n$ is densely defined if and
only if $\escr \subseteq \dz{\slam^n}$.
   \end{enumerate}
   \end{thm}
   \begin{proof}
   (ii) We proceed by induction on $n$. The case of
$n=0$ is obvious. Assume that assertion (ii) holds for
a fixed $n \in \zbb_+$. Take $f$ in
$\dz{\slam^{n+1}}$. It follows from \eqref{roz} that
   \begin{align}  \label{lamtauf+}
\{v \in V^\circ\colon f(\paa(v)) \neq 0, \, \lambda_v
\neq 0\} = \bigsqcup_{u \in V :\, f(u) \neq 0} \{v \in
\dzi{u}\colon \lambda_v \neq 0\}.
   \end{align}
Applying the induction hypothesis to the function
$\slam f$ which clearly belongs to $\dz{\slam^n}$, we
obtain
   \allowdisplaybreaks
   \begin{align} \label{series}
\slam^{n+1} f = \slam^n(\slam f)
&\overset{\eqref{aaa}}= \underset{v\in V:\,\, (\slam
f)(v) \neq 0}{\sumo} (\slam f)(v) \slam^n e_v
   \\ \notag
&\overset{\eqref{lamtauf}}= \underset{v\in V^\circ:\,
f(\paa(v))\neq 0, \, \lambda_v \neq 0}{\sumo}
\lambda_v f(\paa(v))\slam^n e_v
   \\ \notag
&\overset{\eqref{lamtauf+}} = \underset{u\in V:\, f(u)
\neq 0}{\sumo} f(u) \Bigg(\,\underset{v \in \dzi{u}:\,
\lambda_v \neq 0}{\sumo} \lambda_v\slam^n e_v\Bigg)
   \\ \notag
&\overset{\eqref{pow}}= \underset{u\in V:\, f(u) \neq
0}{\sumo} f(u) \Bigg(\,\underset{v \in \dzi{u}:\,
\lambda_v \neq 0}{\sumo} \;\underset{w \in
\dzin{n}{v}}{\sumo} \lambda_v \lambda_{v\mid w} e_w
\Bigg)
   \\  \notag
&\overset{\eqref{recfor2}} = \underset{u\in V:\, f(u)
\neq 0}{\sumo} f(u) \Bigg(\,\underset{v \in \dzi{u}:\,
\lambda_v \neq 0}{\sumo} \; \underset{w \in
\dzin{n}{v}}{\sumo} \lambda_{u\mid w} e_w\Bigg)
   \\  \notag
&\overset{\eqref{dzinn2}}= \underset{u\in V:\, f(u)
\neq 0}{\sumo} f(u) \Bigg(\,\underset{w \in
\dzin{n+1}{u}:\, \lambda_{\paa^n(w)} \neq 0}{\sumo}
\lambda_{u\mid w} e_w\Bigg)
   \\  \notag
&\overset{\eqref{luv}} = \underset{u\in V:\, f(u) \neq
0}{\sumo} f(u) \Bigg(\,\underset{w \in
\dzin{n+1}{u}}{\sumo} \lambda_{u\mid w} e_w\Bigg),
   \end{align}
where the penultimate inequality is valid because the
vectors $\{e_w\}_{w \in V}$ are pairwise orthogonal.
Since the series in \eqref{series} are orthogonal, we
deduce that
   \begin{align} \label{brak2}
\sum_{u\in V} |f(u)|^2 \sum_{v \in \dzin{n+1}{u}}
|\lambda_{u\mid v}|^2 = \|\slam^{n+1} f\|^2,
   \end{align}
and hence that
   \begin{align} \label{ciach}
\text{$\sum_{v \in \dzin{n+1}{u}} |\lambda_{u\mid
v}|^2 < \infty$ for every $u \in V$ such that
$f(u)\neq 0$.}
   \end{align}
As $f$ belongs to $\dz{\slam^n}$, we infer from Lemma
\ref{lem4}\,(i) and the induction hypothesis applied
to $f$ that $\sum_{v \in \dzin{m}{u}} |\lambda_{u\mid
v}|^2 < \infty$ for $m=0, \ldots, n$ and for every $u
\in V$ such that $f(u)\neq 0$. But this, together with
\eqref{ciach} and Lemma \ref{lem4}\,(i), implies that
$e_u \in \dz{\slam^{n+1}}$ for all $u\in V$ such that
$f(u)\neq 0$. Combining Lemma \ref{lem4}\,(ii) with
\eqref{series} and \eqref{brak2}, we obtain
\eqref{aaa} and \eqref{bbb} with $n+1$ in place of
$n$, which completes the induction argument.
Therefore, (ii) holds.

   (i) It follows from (ii) and \eqref{brak} that the
``only if'' part of assertion (i) holds for all $n\in
\zbb_+$. To prove the reverse implication in (i), we
proceed by induction on $n$. The case of $n=0$ is
obvious. Assume that for a fixed $n \in \zbb_+$, the
``if'' part of assertion (i) holds. Let $f\colon V \to
\cbb$ be a function satisfying \eqref{series2} with
$n+1$ in place of $n$. Since $n+1\Ge 1$, this implies
that
   \begin{align*}
\sum_{u\in V} |f(u)|^2 \Big(1 + \sum_{v \in \dzi{u}}
|\lambda_{v}|^2\Big) < \infty,
   \end{align*}
which in view of \cite[Proposition 3.1.3(i)]{j-j-s}
yields $f \in \dz{\slam}$. Note that
   \allowdisplaybreaks
   \begin{align*}
\sum_{u \in V} |(\slam f)(u)|^2 & \Bigg(\sum_{j=0}^n
\, \sum_{v \in \dzin{j}{u}} |\lambda_{u\mid
v}|^2\Bigg)
   \\
& \overset{\eqref{lamtauf}}= \sum_{u \in V^\circ}
|\lambda_u f(\paa(u))|^2 \Bigg(\sum_{j=0}^n \, \sum_{v
\in \dzin{j}{u}} |\lambda_{u\mid v}|^2\Bigg)
   \\
& \overset{\eqref{roz}}= \sum_{x\in V} \sum_{u \in
\dzi{x}} |\lambda_u|^2 |f(x)|^2 \Bigg(\sum_{j=0}^n \,
\sum_{v \in \dzin{j}{u}} |\lambda_{u\mid v}|^2\Bigg)
   \\
&\hspace{1.7ex}= \sum_{x\in V} |f(x)|^2 \sum_{j=0}^n
\sum_{u \in \dzi{x}} \sum_{v \in \dzin{j}{u}}
|\lambda_u \lambda_{u\mid v}|^2
   \\
& \overset{\eqref{recfor2}}= \sum_{x\in V} |f(x)|^2
\sum_{j=0}^n \sum_{u \in \dzi{x}} \sum_{v \in
\dzin{j}{u}} |\lambda_{x\mid v}|^2
   \\
& \overset{\eqref{dzinn2}} = \sum_{x\in V} |f(x)|^2
\sum_{j=0}^n \sum_{v \in \dzin{j+1}{x}}
|\lambda_{x\mid v}|^2 < \infty.
   \end{align*}
Hence, by the induction hypothesis, we see that $\slam
f$ is in $\dz{\slam^n}$. This completes the proof of
(i).

   (iii) Suppose that $\escr\subseteq \dz{\slam^n}$.
Thus, by (i), (ii) and \eqref{brak}, the domain and
the graph norm of $\slam^n$ are given by the following
formulas:
   \begin{gather*}
\dz{\slam^n} = \Big\{f\in \cbb^V \colon \sum_{u\in V}
|f(u)|^2 \Big(\sum_{j=0}^n \, \|\slam^j e_u\|^2\Big) <
\infty\Big\},
   \\
\|f\|^2 + \|\slam^n f\|^2 = \sum_{u\in V} |f(u)|^2
\Big(1 + \|\slam^n e_u\|^2\Big), \quad f \in
\dz{\slam^n}.
   \end{gather*}
Since $\escr$, being the set of all complex functions
on $V$ which vanish off finite sets, is dense in the
weighted $\ell^2$-space on $V$ with weights $\big\{1 +
\|\slam^n e_u\|^2\big\}_{u \in V}$ and $\dz{\slam^n}$
is between these two spaces, we see that $\escr$ is a
core of $\slam^n$.

(iv) Since $\escr$ is dense in $\ell^2(V)$, we see
that the ``if'' part of assertion (iv) is valid.
Suppose that the reverse implication in (iv) does not
hold. Then $\slam^n$ is densely defined and $e_u
\notin \dz{\slam^n}$ for some $u \in V$. Hence, by
(ii), $f(u)=0$ for every $f \in \dz{\slam^n}$. This
and the density of $\dz{\slam^n}$ in $\ell^2(V)$ imply
that $e_u \perp \ell^2(V)$, which is a contradiction.
This completes the proof of Theorem \ref{dn}.
   \end{proof}
Regarding Theorem \ref{dn}, we note that classical
unilateral and bilateral weight\-ed shifts are always
closed, but their higher powers may not be closed.
   \begin{cor} \label{polem}
Let $\slam$ be a weighted shift on a directed tree
$\tcal$ with weights $\lambdab = \{\lambda_v\}_{v\in
V^\circ}$. Then the following conditions are
equivalent{\/\em :}
   \begin{enumerate}
   \item[(i)] $\escr \subseteq \dzn{\slam}$,
   \item[(ii)] $\dzn{\slam}$ is dense in $\ell^2(V)$,
   \item[(iii)] $\slam^n$ is densely defined for
every $n\in \zbb_+$.
   \end{enumerate}
Moreover, if any of the above equivalent conditions
holds, then $\dzn{\slam}$ is a core of $\slam^n$ for
every $n\in\zbb_+$.
   \end{cor}
It is worth pointing out that the equivalence
(ii)$\Leftrightarrow$(iii) which appears in Corollary
\ref{polem} remains true in the class of composition
operators in $L^2$-spaces (cf.\ \cite{b-j-j-s2}).

We conclude this section by proving that a weighted
shift $\slam$ on a directed tree generates Stieltjes
moment sequences if and only if each basic vector
$e_u$, $u \in V$, induces a Stieltjes moment sequence.
   \begin{thm} \label{cinf}
Let $\slam$ be a weighted shift on a directed tree
$\tcal$ with weights $\lambdab = \{\lambda_v\}_{v\in
V^\circ}$. Suppose that $\escr \subseteq \dzn{\slam}$
and $\{\|\slam^n e_u\|^2\}_{n=0}^\infty$ is a
Stieltjes moment sequence for every $u \in V$. Then
$\{\|\slam^n f\|^2\}_{n=0}^\infty$ is a Stieltjes
moment sequence for every $f \in \dzn{\slam}$.
   \end{thm}
   \begin{proof}
Since, by the Stieltjes theorem (cf.\ \cite[Theorem
6.2.5]{ber}), the class of Stieltjes moment sequences
is closed under both, the operation of taking linear
combinations with nonnegative coefficients and the
operation of taking pointwise limits, we can infer
Theorem \ref{cinf} from \eqref{bbb}.
   \end{proof}
   \section{EXAMPLES OF EXOTIC NON-HYPONORMAL OPERATORS}
   \subsection{\label{sec8}General scheme}
In this section we introduce a class of weighted
shifts on an enumerable leafless directed tree with
one branching vertex. Such a directed tree (which is,
roughly speaking, one step more complicated than the
directed trees involved in the definitions of
classical weighted shifts) can be modelled as follows
(cf.\ \cite[(6.2.10)]{j-j-s}). Given $\eta,\kappa \in
\zbb_+ \sqcup \{\infty\}$ with $\eta \Ge 2$, we define
the directed tree $\tcal_{\eta,\kappa} =
(V_{\eta,\kappa}, E_{\eta,\kappa})$ by
   \allowdisplaybreaks
   \begin{align*}
   \begin{aligned}
V_{\eta,\kappa} & = \big\{-k\colon k\in J_\kappa\big\}
\sqcup \{0\} \sqcup \big\{(i,j)\colon i\in J_\eta,\,
j\in \nbb\big\},
   \\
E_{\eta,\kappa} & = E_\kappa \sqcup
\big\{(0,(i,1))\colon i \in J_\eta\big\} \sqcup
\big\{((i,j),(i,j+1))\colon i\in J_\eta,\, j\in
\nbb\big\},
   \\
E_\kappa & = \big\{(-k,-k+1) \colon k\in
J_\kappa\big\},
   \end{aligned}
   \end{align*}
where
   \begin{align*}
J_\iota = \{k \in \nbb\colon k\Le \iota\}, \quad \iota
\in \zbb_+ \sqcup \{\infty\}.
   \end{align*}
Note that $0$ is the only branching vertex of
$\tcal_{\eta,\kappa}$ and $V_{\eta,\kappa}^\circ =
V_{\eta,\kappa} \setminus \{-\kappa\}$.

Let $\gamkap{-\kappa}$ be a system of positive real
numbers such that
   \begin{align} \label{g1}
\gamma_0 & =1
   \\
\gamma_n & = \int_0^\infty x^n \D \nu(x), \quad n \in
\zbb, \, n \Ge -\kappa, \label{g2}
   \end{align}
for some positive Borel measure $\nu$ on $\rbb_+$
(note that if $\kappa > 0$, then \eqref{g2} implies
that $\nu(\{0\})=0$). It follows from \eqref{g2} that
   \begin{align} \label{g3}
\text{$\{\gamma_{n-k}\}_{n=0}^\infty$ is a Stieltjes
moment sequence for every integer $k\Le \kappa$.}
   \end{align}
Suppose that there exists an \Sti-representing measure
$\rho$ of $\gamsj$ such that
   \allowdisplaybreaks
   \begin{align}     \label{jjs1}
&0 < \int_0^\infty \frac{1}{x^n} \D \rho(x) < \infty,
\quad n \in J_{\kappa+1},
   \\  \label{jjs2}
&\card{\supp{\rho}} \Ge
   \begin{cases} \eta & \text{
if } \eta < \infty,
   \\
\aleph_0 & \text{ if } \eta=\infty.
   \end{cases}
   \end{align}
Let $\{\varOmega_i\}_{i=1}^\eta$ be sequence of
pairwise disjoint Borel subsets of $(0,\infty)$ such
that
   \begin{align} \label{mu1}
\rho(\varOmega_i) & > 0, \quad i\in J_\eta,
   \\       \label{mu2}
\bigsqcup_{i\in J_\eta} \varOmega_i &= (0,\infty).
   \end{align}
Since, by \eqref{jjs1}, $0$ is not an atom of $\rho$,
one can deduce from \eqref{jjs2} that such
$\{\varOmega_i\}_{i=1}^\eta$ always exists (see also
Proposition \ref{atom} for the case of
$\card{\supp{\rho}}=\aleph_0$). In view of
\eqref{mu1}, we can define the sequence
$\{\mu_{i,1}\}_{i\in J_\eta}$ of Borel probability
measures on $\rbb_+$~ by
   \begin{align} \label{jjs3}
\mu_{i,1}(\sigma) = \frac{1}{\rho(\varOmega_i)}
\rho(\varOmega_i \cap \sigma), \quad \sigma \in
\borel{\rbb_+}, \, i\in J_\eta,
   \end{align}
and the family $\{\lambda_{i,j}\colon i\in J_\eta, j
\in \nbb\}$ of positive real numbers by
   \begin{align} \label{eq4}
\lambda_{i,j} =
   \begin{cases}
\sqrt{\rho(\varOmega_i)} & \text{ for } j=1,
   \\[2ex]
\sqrt{\cfrac{\int_0^\infty x^{j-1} \D
\mu_{i,1}(x)}{\int_0^\infty x^{j-2} \D \mu_{i,1}(x)}}
& \text{ for } j\Ge 2,
   \end{cases}
\quad i\in J_\eta.
   \end{align}
If $\kappa > 0$, then we define the sequence of
$\{\lambda_{-k}\}_{k=0}^{\kappa-1}$ of positive real
numbers by
   \begin{align}  \label{lambda-k}
\lambda_{-k} & =
\sqrt{\frac{\gamma_{-k}}{\gamma_{-(k+1)}}}, \quad k
\in \zbb_+,\, 0 \Le k < \kappa.
   \end{align}
Let $\slam$ be a weighted shift on the directed tree
$\tcal_{\eta,\kappa}$ with weights
$\lambdab=\{\lambda_v\}_{v\in V_{\eta,\kappa}^\circ}$
defined by \eqref{eq4} and \eqref{lambda-k} (we adhere
to notation $\lambda_{i,j}$ instead of a more formal
expression $\lambda_{(i,j)}$). The reader should be
aware of the fact that the operator $\slam$ just
constructed depends not only on $\gamkap{-\kappa}$ and
$\rho$, but also on the partition
$\{\varOmega_i\}_{i=1}^\eta$ of $(0,\infty)$. Now we
can prove some crucial properties of $\slam$.
   \begin{thm} \label{import}
Let $\gamkap{-\kappa}$, $\rho$, $\{\varOmega_i\}_{i\in
J_\eta}$, $\{\mu_{i,1}\}_{i\in J_\eta}$,
$\lambdab=\{\lambda_v\}_{v\in V_{\eta,\kappa}^\circ}$
and $\slam$ be as above. Then the following assertions
hold.
   \begin{enumerate}
   \item[(i)] $\EuScript{E}_{V_{\eta,\kappa}} \subseteq
\dzn{\slam}$.
   \item[(ii)] $\{\|\slam^n f\|^2\}_{n=0}^\infty$ is
a Stieltjes moment sequence for every $f \in
\dzn{\slam}$.
   \item[(iii)] $\slam$ is paranormal.
   \item[(iv)] The consistency condition
\eqref{con2} holds at $u=0$ if and only if
   \begin{align} \label{fura}
\int_0^\infty \frac 1 x\, \D \rho(x) \Le 1.
   \end{align}
   \item[(v)] $\slam$ is hyponormal if and only if
    \begin{gather} \label{hypon2}
\sum_{i\in J_\eta} \frac {\lambda_{i,1}^2}{\|\slam
e_{i,1}\|^2} \Le 1.
    \end{gather}
   \item[(vi)] The following inequality holds
   \begin{align} \label{hihi}
\sum_{i\in J_\eta} \frac {\lambda_{i,1}^2}{\|\slam
e_{i,1}\|^2} \Le \int_0^\infty \frac{1}{x} \D\rho(x).
   \end{align}
   \item[(vii)] The inequality in \eqref{hihi}
turns into equality if and only if for every $i\in
J_\eta$, there exists $q_i \in \varOmega_i$ such that
   \begin{align} \label{fura2}
\rho(\sigma \cap \varOmega_i) = \rho(\varOmega_i)
\cdot \delta_{q_i}(\sigma), \quad \sigma \in
\borel{\rbb_+}, \, i \in J_\eta.
   \end{align}
   \item[(viii)] If $\gamsj$ is \Sti-determinate, then
$\slam$ is subnormal and \eqref{fura} holds.
   \end{enumerate}
   \end{thm}
   \begin{proof}
   We prove (i) and (ii) simultaneously. It follows
from Lemma \ref{lem4}(i) that $e_{i,j} \in
\dzn{\slam}$ for all $(i,j) \in \des{0} \setminus
\{0\}$ (we abbreviate $e_{(i,j)}$ to $e_{i,j}$), and
   \begin{align}   \label{jjs4}
\|\slam^n e_{i,j}\|^2 & \overset{\eqref{eq4}}=
\int_0^\infty x^n \D \mu_{i,j}(x), \quad n \in
\zbb_+, \, (i,j) \in \des{0} \setminus \{0\},
   \end{align}
where
   \begin{align*}
\mu_{i,j}(\sigma) & \hspace{1ex}=
\frac{1}{\int_0^\infty x^{j-1} \D \mu_{i,1}(x)}
\int_{\sigma} x^{j-1} \D \mu_{i,1}(x), \quad
\sigma \in \borel{\rbb_+}, \, (i,j) \in \des{0}
\setminus \{0\}.
   \end{align*}
Noting that
   \begin{multline*}
\sum_{v \in \dzin{n}{0}} \lambda_{0\mid v}^2
\hspace{1ex}= \sum_{i\in J_\eta} \lambda_{0\mid
(i,n)}^2 \overset{\eqref{luv}} = \sum_{i \in J_\eta}
\prod_{j=1}^n \lambda_{i,j}^2
   \\
\overset{\eqref{jjs3}\& \eqref{eq4}}= \sum_{i \in
J_\eta} \int_{\varOmega_i} x^{n-1} \D\rho(x)
\overset{\eqref{mu2}} = \int_0^\infty x^{n-1}
\D\rho(x) < \infty, \quad n \Ge 1,
   \end{multline*}
and applying Lemma \ref{lem4}, we deduce that $e_0 \in
\dzn{\slam}$ and
   \begin{align}      \label{jjs6}
\|\slam^{n+1} e_0\|^2 = \int_0^\infty x^{n} \D\rho(x),
\quad n \in \zbb_+.
   \end{align}
As $\rho$ is an \Sti-representing measure of $\gamsj$,
we infer from \eqref{g1} and \eqref{jjs6} that
   \begin{align}   \label{jjs5}
\|\slam^{n} e_0\|^2 = \gamma_n, \quad n \in \zbb_+.
   \end{align}
Now combining \eqref{jjs4} with \eqref{jjs5}, we
conclude that $\{\|\slam^n e_u\|^2\}_{n=0}^\infty$ is
a Stieltjes moment sequence for every $u \in \des{0}$.

Consider now the case of $\kappa > 0$. By using Lemma
\ref{lem4} and the fact that $e_0 \in \dzn{\slam}$, we
deduce that $e_{-k} \in \dzn{\slam}$ for every $k \in
J_\kappa$, which means that (i) holds. Now we show
that
   \begin{align} \label{g5}
\|\slam^n e_{-k}\|^2 =
\frac{\gamma_{n-k}}{\gamma_{-k}}, \quad n \in
\zbb_+, \, k \in J_\kappa.
   \end{align}
Indeed, if $k \in J_\kappa$, then
   \begin{align} \label{+++}
\|\slam^n e_{-k}\|^2 \overset{\eqref{brak}}=
\prod_{j=k-n}^{k-1} \lambda_{-j}^2
\overset{\eqref{lambda-k}}= \prod_{j=k-n}^{k-1}
\frac{\gamma_{-j}}{\gamma_{-(j+1)}} =
\frac{\gamma_{n-k}}{\gamma_{-k}}, \quad n \in J_k,
   \end{align}
which, in view of \eqref{jjs5} and \eqref{g1}, yields
   \begin{align} \label{4+}
\|\slam^n e_{-k}\|^2 = \prod_{j=0}^{k-1}
\lambda_{-j}^2 \|\slam^{n-k} e_{0}\|^2 =
\frac{\gamma_0}{\gamma_{-k}} \gamma_{n-k} =
\frac{\gamma_{n-k}}{\gamma_{-k}}, \quad n \in \zbb, \,
n > k.
   \end{align}
Combining \eqref{+++} with \eqref{4+}, we obtain
\eqref{g5}. It follows from \eqref{g3} and \eqref{g5}
that the sequence $\{\|\slam^n
e_{-k}\|^2\}_{n=0}^\infty$ is a Stieltjes moment
sequence for every $k\in J_\kappa$. Together with
\eqref{jjs4} and \eqref{jjs5}, this implies that
$\{\|\slam^n e_u\|^2\}_{n=0}^\infty$ is a Stieltjes
moment sequence for every $u \in V_{\eta,\kappa}$.
Thus, by Theorem \ref{cinf}, assertion (ii) is proved.

   (iii) Fix $h \in \dzn{\slam}$. Then, by (ii), there
exists a positive Borel measure $\mu_h$ on $\rbb_+$
such that $\|\slam^n h\|^2 =\int_0^\infty x^n \D
\mu_h(x)$ for all $n\in \zbb_+$. By the Cauchy-Schwarz
inequality, we have
   \begin{align}  \label{para}
   \begin{aligned}
\|\slam h\|^2 & = \int_0^\infty x^0 x^1 \D \mu_h(x)
   \\
&\Le \Big(\int_0^\infty x^0 \D \mu_h(x)\Big)^{\frac
12} \Big(\int_0^\infty x^2 \D \mu_h(x)\Big)^{\frac 12}
= \|h\| \|\slam^2 h\|.
   \end{aligned}
   \end{align}
Take $f \in \dz{\slam^2}$. It follows from (i) and
Corollary \ref{polem} that there exists a sequence
$\{h_n\}_{n=1}^\infty \subseteq \dzn{\slam}$ such that
$h_n \to f$ and $\slam^2 h_n \to \slam^2 f$ as $n\to
\infty$. This and \eqref{para} yield
   \begin{align*}
\|\slam h_m - \slam h_n\|^2 \Le \|h_m - h_n\|
\|\slam^2 h_m - \slam^2 h_n\|, \quad m,n \in \nbb,
   \end{align*}
which, by the completeness of $\hh$, implies that the
sequence $\{\slam h_n\}_{n=1}^\infty$ is convergent in
$\hh$. Since $\slam$ is closed (cf.\ Proposition
\ref{bas}(i)), we deduce that $\slam h_n \to \slam f$
as $n \to \infty$. Hence, by passage to the limit in
the inequality $\|\slam h_n\|^2 \Le \|h_n\| \|\slam^2
h_n\|$ (see \eqref{para}), we obtain $\|\slam f\|^2
\Le \|f\| \|\slam^2 f\|$. This shows that $\slam$ is
paranormal.

(iv) Since $\rho(\{0\})=0$, we obtain
   \begin{align*}
\sum_{v \in \dzi{0}} \lambda_v^2 \int_0^\infty \frac 1
x\, \D \mu_v(x) & = \sum_{i \in J_\eta}
\lambda_{i,1}^2 \int_0^\infty \frac 1 x\, \D
\mu_{i,1}(x)
   \\
& \hspace{-5ex} \overset{\eqref{jjs3}\&\eqref{eq4}}=
\sum_{i \in J_\eta} \int_{\varOmega_i} \frac 1 x\, \D
\rho(x) \overset{\eqref{mu2}}= \int_0^\infty \frac 1
x\, \D \rho(x),
   \end{align*}
which yields (iv).

   (v) Inequality \eqref{hypon}, written for
$u=(i,j) \in \des{0} \setminus \{0\}$, takes the
form $\lambda_{i,j+1} \Le \lambda_{i,j+2}$, which
in view of \eqref{eq4} is equivalent to
   \begin{align} \label{argu}
   \begin{aligned}
\Big(\int_0^\infty x^{j} \D \mu_{i,1}(x)\Big)^2 &
= \Big(\int_0^\infty \sqrt{x^{j-1}}
\sqrt{x^{j+1}} \D \mu_{i,1}(x)\Big)^2
   \\
& \Le \int_0^\infty x^{j-1} \D \mu_{i,1}(x)
\int_0^\infty x^{j+1} \D \mu_{i,1}(x).
   \end{aligned}
   \end{align}
Since the latter is always true due to the
Cauchy-Schwarz inequality, we see that \eqref{hypon}
is valid for all $u \in \des{0} \setminus \{0\}$.
Clearly, inequality \eqref{hypon} is valid for $u=0$
if and only if \eqref{hypon2} holds. Finally, if
$\kappa > 0$ and $k \in J_\kappa$, then using
\eqref{g1} and \eqref{g2} and arguing as in
\eqref{argu}, we verify that $\gamma_{-(k-1)}^2 \Le
\gamma_{-k} \gamma_{-(k-2)}$ for any integer $k$ such
that $2 \Le k \Le \kappa$, and that $\gamma_0 =
\gamma_0^2 \Le \gamma_{-1} \gamma_1$. Hence, by
\eqref{lambda-k} and \eqref{jjs5} applied to $n=1$, we
conclude that inequality \eqref{hypon} is valid for
$u=-k$ whenever $k \in J_\kappa$. Applying Theorem
\ref{hyp} yields (v).

   (vi) It follows from \eqref{mu2} and the
Cauchy-Schwarz inequality that
   \begin{align}    \label{3hihi}
\rho(\varOmega_i)^2 = \Big(\int_{\varOmega_i}
\frac{1}{\sqrt{x}} \sqrt{x} \D \rho(x)\Big)^2 \Le
\int_{\varOmega_i} \frac{1}{x} \D\rho(x) \cdot
\int_{\varOmega_i} x \D\rho(x), \quad i \in J_\eta,
   \end{align}
which together with \eqref{mu1} implies that
   \begin{align}  \label{nier1}
\frac{\rho(\varOmega_i)^2}{\int_{\varOmega_i} x
\D\rho(x)} \Le \int_{\varOmega_i} \frac{1}{x}
\D\rho(x), \quad i \in J_\eta.
   \end{align}
Therefore, by \eqref{jjs3}, \eqref{eq4} and
\eqref{jjs4}, we have (recall that $\rho(\{0\})=0$)
   \begin{align} \label{2hihi}
   \begin{aligned} \sum_{i\in J_\eta}
\frac {\lambda_{i,1}^2}{\|\slam e_{i,1}\|^2} & =
\sum_{i\in J_\eta} \frac
{\rho(\varOmega_i)^2}{\int_{\varOmega_i} x \D \rho(x)}
   \\
& \hspace{-2.2ex}\overset{\eqref{nier1}} \Le
\sum_{i\in J_\eta} \int_{\varOmega_i} \frac{1}{x}
\D\rho(x) \overset{\eqref{mu2}}= \int_0^\infty
\frac{1}{x} \D\rho(x),
   \end{aligned}
   \end{align}
which gives (vi).

(vii) If we have equality in \eqref{hihi}, then one
can deduce from \eqref{nier1} and \eqref{2hihi} that
the inequality in \eqref{nier1} turns into equality
for every $i \in J_\eta$. The latter is equivalent to
the fact that the Cauchy-Schwarz inequality
\eqref{3hihi} becomes an equality for every $i \in
J_\eta$. Since this is possible if and only if the
functions $\frac{1}{\sqrt{x}}$ and $\sqrt{x}$ are
linearly dependent as vectors in
$L^2(\varOmega_i,\borel{\varOmega_i},\rho)$ for every
$i \in J_\eta$, we conclude that \eqref{fura2} holds
for some sequence $\{q_i\}_{i\in J_\eta}$ such that
$q_i \in \varOmega_i$ for all $i\in J_\eta$. The
reverse implication is obvious.

   (viii) Since the Stieltjes moment sequence $\gamsj$
is \Sti-determinate, we infer from \eqref{jjs5} that
$\{\|\slam^{n+1} e_0\|^2\}_{n=0}^\infty$ is an
\Sti-determinate Stieltjes moment sequence. This fact
together with (i) and (ii) implies that the weighted
shift $\slam$ (which has nonzero weights) satisfies
all the assumptions of \cite[Corollary
6.2.5]{b-j-j-s}. Hence, by this corollary, $\slam$ is
subnormal and it satisfies the consistency condition
\eqref{con2} at $u=0$. Applying (iv) completes the
proof.
   \end{proof}
Note that, in virtue of Theorem \ref{import}, the
validity of the consistency condition \eqref{con2} at
$u=0$ implies the hyponormality of $\slam$.

Regarding Theorem \ref{import}\,(vii) it is worth
mentioning that if $\card{\supp{\rho}} = \aleph_0$,
then we can always find a Borel partition
$\{\varOmega_i\}_{i=1}^\infty$ of $(0,\infty)$
satisfying \eqref{mu1} and \eqref{fura2} with
$\eta=\infty$.
   \begin{pro}\label{atom}
Let $\rho$ be a finite positive Borel measure on
$\rbb_+$ such that $\card{\supp{\rho}}=\aleph_0$. Then
there exist a Borel partition
$\{\varOmega_i\}_{i=1}^\infty$ of $(0,\infty)$ and a
sequence $\{q_i\}_{i=1}^\infty \subseteq (0,\infty)$
such that $\rho(\varOmega_i) > 0$, $q_i \in
\varOmega_i$ and $\rho(\sigma \cap \varOmega_i) =
\rho(\varOmega_i) \cdot \delta_{q_i}(\sigma)$ for all
$\sigma \in \borel{\rbb_+}$ and $i \in \nbb$.
   \end{pro}
   \begin{proof}
Clearly $\supp{\rho} = A \sqcup B$, where $A:= \{x \in
\supp{\rho}\colon \rho(\{x\}) > 0\}$ and $B :=
\supp{\rho} \setminus A$. Since
$\card{\supp{\rho}}=\aleph_0$, we deduce that $B
\subseteq A^\prime$, where $A^\prime$ is the set of
all accumulation points of $A$ in $\rbb_+$. This and
the equality $\card{\supp{\rho}}=\aleph_0$ imply that
$\card{A}=\aleph_0$. Hence there exists a sequence
$\{q_i\}_{i=1}^\infty$ of distinct positive real
numbers such that $A \setminus \{0\} = \{q_1, q_2,
q_3, \ldots\}$. Set
   \begin{align*}
\varOmega_i=
   \begin{cases}
   \big((0,\infty) \setminus \supp{\rho}\big) \sqcup
\{q_1\} & \text{ if } i=1,
   \\
   \big(B \setminus \{0\}\big) \sqcup \{q_2\} & \text{
if } i=2,
   \\
   \{q_i\} & \text{ if } i \Ge 3.
   \end{cases}
   \end{align*}
It is a simple matter to verify that
$\{\varOmega_i\}_{i=1}^\infty$ is the required Borel
partition of $(0,\infty)$, which completes the proof.
   \end{proof}
   \subsection{\label{sec9}The main example}
   The following example was announced in the title of
this paper.
   \begin{exa} \label{exnindeterm}
Fix $\kappa \in \zbb_+ \sqcup \{\infty\}$. Let
$\{\gamma_n\}_{n=-\kappa}^\infty$, $\nu$ and $\rho$ be
as in Example \ref{e1}, i.e.,
$\{\gamma_n\}_{n=-\kappa}^\infty$ is a system of
positive real numbers and $\nu, \rho$ are positive
Borel measures on $\rbb_+$ satisfying the conditions
(i) to (iv) of this example. From (iii) and (iv) we
infer that $\card{\supp{\rho}} = \aleph_0$. Hence the
triplet $(\{\gamma_n\}_{n=-\kappa}^\infty, \nu, \rho)$
satisfies the conditions \eqref{g1}, \eqref{g2},
\eqref{jjs1} and \eqref{jjs2} with $\eta=\infty$. It
follows from Proposition \ref{atom} that there exist a
Borel partition $\{\varOmega_i\}_{i=1}^\infty$ of
$(0,\infty)$ and a sequence $\{q_i\}_{i=1}^\infty
\subseteq (0,\infty)$ which satisfy \eqref{mu1} and
\eqref{fura2} with $\eta=\infty$ (note that if
$\rho=\tilde \rho_a$ for some $a \in (0,\infty)$,
where $\tilde \rho_a$ is as in Example \ref{e1}, then
we may simply consider the sequence
$\{q_i\}_{i=1}^\infty := \{a, aq, aq^{-1}, aq^2, a
q^{-2}, \ldots\}$). Let $\slam$ be a weighted shift on
the directed tree $\tcal_{\infty,\kappa}$ with weights
$\lambdab=\{\lambda_v\}_{v\in
V_{\infty,\kappa}^\circ}$ defined by \eqref{jjs3},
\eqref{eq4} and \eqref{lambda-k} with $\eta=\infty$.
By \eqref{ineq1} and assertions (v), (vi) and (vii) of
Theorem \ref{import}, the operator $\slam$ is not
hyponormal. In turn, assertions (i), (ii) and (iii) of
Theorem \ref{import} imply that $\slam$ is a
paranormal operator which generates Stieltjes moment
sequences; moreover, by Corollary \ref{polem},
$\dzn{\slam}$ is a core of $\slam^n$ for every
$n\in\zbb_+$. In view of \eqref{ineq1} and assertion
(iv) of Theorem \ref{import}, the weighted shift
$\slam$ does not satisfy the consistency condition
\eqref{con2} at $u=0$. Since
$\EuScript{E}_{V_{\infty,\kappa}} \subseteq
\dzn{\slam}$ and $\slam$ is not subnormal, we deduce
from \cite[Theorem 5.2.1]{b-j-j-s} that the weighted
shift $\slam$ has no consistent system of measures (in
the sense of \cite{b-j-j-s}). Finally, by making an
appropriate choice of the triplet
$(\{\gamma_n\}_{n=-\kappa}^\infty, \nu, \rho)$, we can
guarantee that $\{\|\slam^{n+1}e_0\|^2\}_{n=0}^\infty$
is \Sti-indeterminate, while
$\{\|\slam^{n}e_0\|^2\}_{n=0}^\infty$ is either
\Ham-determinate or \Sti-indeterminate according to
our needs (cf.\ Example \ref{e1}).
   \end{exa}
The directed tree $\tcal_{\eta,\kappa}$ can also be
used to construct examples of unbounded subnormal
weighted shifts $\slam$ for which
$\{\|\slam^{n}e_0\|^2\}_{n=0}^\infty$ and
$\{\|\slam^{n+1}e_0\|^2\}_{n=0}^\infty$ are
\Ham-determinate Stieltjes moment sequences.
   \begin{exa}
The following example is an adaptation of
\cite[Example 7.1]{js-jbs} to our needs. Set $c =
\sum_{j=2}^\infty j2^{-j} + \sum_{j=2}^\infty
j^{-1}\E^{-j^2}$. It is easily seen that the two-sided
sequence $\{\gamma_n\}_{n=-\infty}^\infty$ given by
   \begin{align*}
\gamma_n = c^{-1} \Big(\sum_{j=2}^\infty \frac{1}{2^j
j^{n-1}} + \sum_{j=2}^\infty
\frac{j^{n-1}}{\E^{j^2}}\Big), \quad n \in \zbb,
   \end{align*}
is well-defined, $\gamma_0=1$ and
   \begin{align} \label{2srho}
\gamma_n=\int_0^\infty x^n \D\nu(x), \quad n \in \zbb,
   \end{align}
where $\nu := c^{-1} \big(\sum_{j=2}^\infty j 2^{-j}
\delta_{\frac{1}{j}} + \sum_{j=2}^\infty j^{-1}
\E^{-j^2} \delta_{j}\big)$. Note also that
   \begin{align}  \label{3srho}
\supp{\nu} = \{0\} \cup \Big\{\ldots, \frac{1}{4},
\frac{1}{3}, \frac{1}{2}\Big\} \cup \{2,3,4, \ldots\}.
   \end{align}
It was proved in \cite[Example 7.1]{js-jbs} that
$\gamma_{2n+1} \Le 4 c^{-1} \, n^n$ for all integers
$n \Ge 4$. This implies that $\gamma_{2n} \Le 5c^{-1}
\, n^n$ for all integers $n \Ge 4$. Hence, by
Carleman's criterion (see e.g., \cite[Corollary
4.5]{sim}), the Stieltjes moment sequences
$\{\gamma_{n}\}_{n=0}^\infty$ and
$\{\gamma_{n+1}\}_{n=0}^\infty$ are \Ham-determinate.
In view of \eqref{2srho}, the positive Borel measure
$$\rho := c^{-1} \big(\sum_{j=2}^\infty 2^{-j}
\delta_{\frac{1}{j}} + \sum_{j=2}^\infty \E^{-j^2}
\delta_{j}\big)$$ is a unique representing measure of
$\{\gamma_{n+1}\}_{n=0}^\infty$. Putting all these
together, we conclude that for every $\eta \in \{2,3,
\dots\} \sqcup \{\infty\}$ and for every $\kappa \in
\zbb_+ \cup \{\infty\}$, the system
$\{\gamma_n\}_{n=-\kappa}^\infty$ and the measures
$\nu$ and $\rho$ satisfy \eqref{g1}, \eqref{g2},
\eqref{jjs1} and \eqref{jjs2}. Take any Borel
partition $\{\varOmega_i\}_{i=1}^\eta$ of $(0,\infty)$
which satisfies \eqref{mu1}. Let $\slam$ be the
weighted shift on the directed tree
$\tcal_{\eta,\kappa}$ with weights $\lambdab =
\{\lambda_v\}_{v \in V_{\eta,\kappa}^\circ}$ defined
by \eqref{jjs3}, \eqref{eq4} and \eqref{lambda-k}.
Then, by assertions (iv) and (viii) of Theorem
\ref{import}, the operator $\slam$ is subnormal and it
satisfies the consistency condition \eqref{con2} at
$u=0$ (in fact, $1 = \gamma_0 = \int_0^\infty \frac 1
x\, \D \rho(x)$). Moreover, by \eqref{jjs5}, the
Stieltjes moment sequences
$\{\|\slam^{n}e_0\|^2\}_{n=0}^\infty$ and
$\{\|\slam^{n+1}e_0\|^2\}_{n=0}^\infty$ are
\Ham-determinate. Note that the operator $\slam$ is
unbounded. Indeed, otherwise by \cite[Notation 6.1.9
and Theorem 6.1.3]{j-j-s}, a unique \Ham-representing
measure of $\{\|\slam^n e_0\|^2\}_{n=0}^\infty$ is
compactly supported. This fact, together with
\eqref{jjs5} and \eqref{2srho}, contradicts
\eqref{3srho}.
   \end{exa}
   \subsection{The case of composition operators}
It turns out that Example \ref{exnindeterm} can be
realized as a composition operator in an $L^2$-space.
Before proving this, we show that a great deal of
weighted shifts on directed trees can be identified
with composition operators in $L^2$-spaces.
   \begin{lem} \label{wsodt2comp}
Let $\slam$ be a weighted shift on a rootless directed
tree $\tcal=(V,E)$ with positive weights
$\lambdab=\{\lambda_v\}_{v \in V^\circ}$. Suppose
$\card{V} = \aleph_0$. Then $\slam$ is unitarily
equivalent to a composition operator $C$ in an
$L^2$-space over a $\sigma$-finite measure space.
Moreover, if the directed tree $\tcal$ is leafless,
then $C$ can be made injective.
   \end{lem}
   \begin{proof}
We begin by proving that for any $(w, \beta) \in V
\times (0,\infty)$ there exists a function
$\alpha\colon \des{w} \to (0,\infty)$ such that
   \begin{align} \label{bdes}
\alpha(w) = \beta \text{ and } \alpha(v) = \lambda_v^2
\alpha(u) \text{ for all } v \in \dzi{u} \text{ and }u
\in \des{w}.
   \end{align}
Indeed, since $\des{w}= \bigsqcup_{n=0}^\infty \dzin n
{w}$ (cf.\ \cite[(2.1.10)]{j-j-s}), we can proceed by
induction. For the base step of the induction, set
$\alpha(v) = \lambda_v^2 \beta$ for $v \in \dzi{w}$.
Fix $n\Ge 1$, and assume that we already have a
function $\alpha \colon \bigsqcup_{j=0}^n \dzin {j}{w}
\to (0,\infty)$ such that $\alpha(w)=\beta$ and
$\alpha(v) = \lambda_v^2 \alpha(u)$ for all $v \in
\dzi{u}$ and $u \in \bigsqcup_{j=0}^{n-1} \dzin
{j}{w}$. Since $\dzin{n+1}{w} = \bigsqcup_{u \in
\dzin{n}{w}} \dzi{u}$ (cf.\ \cite[(6.1.3)]{j-j-s}), we
can extend the function $\alpha$ to
$\bigsqcup_{j=0}^{n+1} \dzin {j}{w}$ by setting
$\alpha(v) = \lambda_v^2 \alpha(u)$ for all $v \in
\dzi{u}$ and $u \in \dzin {n}{w}$. Therefore the
induction step is valid, and so our claim is proved.

Fix $z \in V$. Let $\alpha_0\colon \des{z} \to
(0,\infty)$ be a function satisfying \eqref{bdes} with
$\alpha=\alpha_0$, $w=z$ and $\beta=1$. By
\cite[(3.4)]{j-j-s1}, we have
   \begin{align}   \label{des-des}
\des{\pa{z}} \setminus \des{z} = \{\pa{z}\} \sqcup
\bigsqcup_{w \in \dzi{\pa{z}} \setminus \{z\}}
\des{w}.
   \end{align}
Using \eqref{des-des}, we will extend the function
$\alpha_0$ to a function $\alpha_1\colon \des{\pa{z}}
\to (0,\infty)$ which satisfies \eqref{bdes} with
$\alpha=\alpha_1$, $w=\pa{z}$ and $\beta =
1/\lambda_z^2$. By the preceding paragraph, for every
$w \in \dzi{\pa{z}} \setminus \{z\}$ there exists a
function $\alpha_{1,w}\colon \des{w} \to (0,\infty)$
satisfying \eqref{bdes} with $\alpha=\alpha_{1,w}$ and
$\beta=\lambda_w^2/\lambda_z^2$. Set
$\alpha_1(\pa{z})=1/\lambda_z^2$ and $\alpha_1(v) =
\alpha_{1,w}(v)$ for $v \in \des{w}$ and $w \in
\dzi{\pa{z}} \setminus \{z\}$. Then, by
\eqref{des-des}, the function $\alpha_1$ is
well-defined and it satisfies our requirements. Using
the decomposition $V = \bigcup_{k=0}^\infty
\des{\paa^k(z)}$ (cf.\ \cite[Proposition
2.1.6]{j-j-s}) and induction with $\alpha(\paa^k(z)) =
\prod_{j=0}^{k-1} \lambda_{\paa^j(z)}^{-2}$ for $k\Ge
1$, we get a function $\alpha\colon V \to (0,\infty)$
such that
   \begin{align}   \label{lamfrac}
\alpha(v) = \lambda_v^2 \alpha(u), \quad v \in
\dzi{u}, \, u \in V.
   \end{align}

Define a measure space $(V, \varSigma, \mu)$ by
$\varSigma=2^V$ and $\mu(\{u\})=\alpha(u)$ for every
$u \in V$. Since $\card{V} = \aleph_0$, the measure
$\mu$ is $\sigma$-finite. Let $\phi\colon V \to V$ be
a transformation given by $\phi(u) = \pa{u}$ for all
$u \in V$ ($\phi$ is well-defined because $\tcal$ is
rootless) and let $C$ be a composition operator in
$L^2(\mu)$ defined by
   \begin{align*}
\dz{C} = \{f \in L^2(\mu) \colon f \circ \phi \in
L^2(\mu)\} \text{ and } Cf = f \circ \phi \text{ for }
f \in \dz{C}.
   \end{align*}
If the directed tree $\tcal$ is leafless, then the
transformation $\phi$ is surjective, and thus the
operator $C$ is injective. It is clear that the
operator $C$ is closed\footnote{\;In fact, composition
operators in $L^2$-spaces are always closed (see
\cite{b-j-j-s2}, see also \cite[Lemma 6.2]{ca-hor} for
the case of densely defined composition operators).}.
Now we define the mapping $U \colon \ell^2(V) \to
L^2(\mu)$ by
   \begin{align} \label{Unit}
(Uf)(u) = \frac{f(u)}{\sqrt{\alpha(u)}}, \quad u \in
V, \, f \in \ell^2(V).
   \end{align}
It is easily seen that $U$ is a well-defined unitary
isomorphism such that
   \begin{align}  \label{UslC}
   \begin{aligned}
\big((Uf) \circ \phi\big)(v) &\overset{\eqref{Unit}} =
\frac{f(\phi(v))}{\sqrt{\alpha(\phi(v))}}
   \\
&\overset{\eqref{lamfrac}}= \lambda_v
\frac{f(\pa{v})}{\sqrt{\alpha(v)}}
\overset{\eqref{lamtauf}} = \frac{(\varLambda_\tcal
f)(v)}{\sqrt{\alpha(v)}}, \quad v \in V, \, f \in
\ell^2(V).
   \end{aligned}
   \end{align}
Hence if $f \in \dz{\slam}$, then $\big((Uf) \circ
\phi\big)(v) = (U\slam f)(v)$ for every $v \in V$,
which implies that $Uf \in \dz{C}$ and $CU f = U \slam
f$. This shows that $U\slam \subseteq CU$. In turn, if
$f \in \ell^2(V)$ and $Uf \in \dz{C}$, then by
\eqref{UslC} the function $g\colon V \to \cbb$ given
by
   \begin{align*}
g(v)=\frac{(\varLambda_\tcal f)(v)}{\sqrt{\alpha(v)}},
\quad v \in V,
   \end{align*}
belongs to $L^2(\mu)$. It follows from \eqref{Unit}
that $(U^{-1}g)(v) = (\varLambda_\tcal f)(v)$ for
every $v \in V$, which means that $f \in \dz{\slam}$.
Putting all these together, we conclude that $U \slam
= CU$, or equivalently that $\slam = U^* C U$. This
completes the proof.
   \end{proof}
   \begin{rem}
A close inspection of the proof of Lemma
\ref{wsodt2comp} reveals that if functions $\alpha,
\alpha^\prime\colon V \to (0,\infty)$ satisfy
\eqref{lamfrac}, then there exists $t \in (0,\infty)$
such that $\alpha^\prime(v) = t \alpha(v)$ for all
$v\in V$. If we drop the assumption ``$\card{V} =
\aleph_0$'' in Lemma \ref{wsodt2comp}, then the
composition operator $C$ constructed in its proof acts
in an $L^2$-space over a measure space which is not
necessarily $\sigma$-finite. It follows from
\cite[Proposition 3.1.10]{j-j-s} that if there exists
a densely defined weighted shift on a directed tree
$\tcal$ with nonzero weights, then $\card{V} \Le
\aleph_0$.
   \end{rem}
   The following surprising fact follows directly from
Example \ref{exnindeterm} with $\kappa = \infty$ and
Lemma \ref{wsodt2comp}.
   \begin{thm} \label{compns}
There exists a non-hyponormal composition operator $C$
in an $L^2$-space over a $\sigma$-finite measure space
which is injective, paranormal and which generates
Stieltjes moment sequences. Moreover, $C$ has the
property that $\dzn{C}$ is a core of $C^n$ for every
$n\in\zbb_+$.
   \end{thm}
It is worth pointing out that every composition
operator $C$ in an $L^2$-space over a $\sigma$-finite
measure space which generates Stieltjes moment
sequences has the property that $\dzn{C}$ is a core of
$C^n$ for every $n\in\zbb_+$ (cf.\ \cite{b-j-j-s2}).
   \section{APPENDIX}
%\numberwithin{thm}{section} \setcounter{thm}{0}
%\numberwithin{equation}{section}
%\setcounter{equation}{0}
\subsection{}
As announced in the Introduction, the independence
assertion of Barry Simon's theorem which parameterizes
von Neumann extensions of a closed real symmetric
operator with deficiency indices $(1,1)$ is false (see
Propositions \ref{BS5}, \ref{BS6} and \ref{F2}). For
the reader's convenience, we state the Simon theorem
without typos that appeared in the original version.
We have also added a missing assumption that $\varphi
\neq 0$.
   \begin{caut}
   The reader should be aware of the fact that the
inner products considered in Simon's paper \cite{sim}
are linear in the second factor and anti-linear in the
first. From now on we follow his convention.
   \end{caut}
   \begin{thm}[\mbox{\cite[Theorem 2.6]{sim}}]
\label{BSth} Suppose that $A$ is a closed symmetric
operator so that there exists a complex conjugation
under which $A$ is real. Suppose that $d_+=1$ and that
$\ker(A)=\{0\}$, $\dim\ker(A^*) = 1$. Pick $\varphi
\in \ker(A^*) \setminus \{0\}$, $C \varphi = \varphi$,
and $\eta \in \dz{A^*}$, not in $\dz{A} + \ker(A^*)$.
Then $\is{\varphi}{A^*\eta} \neq 0$ and $\psi =
\{\eta- [\is{\eta}{A^*\eta}/ \is{\varphi}{A^*\eta}]
\varphi\} / \is{\varphi}{A^*\eta}$ are such that in
$\varphi$, $\psi$ basis, $\is{\cdot}{A^* \cdot}$ has
the form
   \begin{align} \tag{2.7}
\is{\cdot}{A^* \cdot} = \left(\begin{matrix} 0 & 1
\\ 0 & 0\end{matrix}\right).
   \end{align}
The self-adjoint extensions, $B_t$, can be labelled by
a real number or $\infty$ where
   \begin{align*}
\dz{B_t} & = \dz{A} + \{\alpha (t\varphi + \psi) \,|
\, \alpha \in \cbb\} & t \in \rbb \hspace{1.5ex}&
   \\
&= \dz{A} + \{\alpha \varphi \;|\, \alpha \in \cbb\} &
t=\infty. &
   \end{align*}
The operators $B_t$ are independent of which real
$\psi$ in $\dz{A^*} \setminus \dz{A}$ is chosen so
that {\em (2.7)} holds.
   \end{thm}
\subsection{}
Let $C$ be a complex conjugation on a complex Hilbert
space $\hh$ (i.e., $C$ is an anti-linear map from
$\hh$ to $\hh$ such that $C(Cf) = f$ and
$\is{Cf}{Cg}=\is{g}{f}$ for all $f,g \in \hh$). We say
that a vector $f$ in $\hh$ is {\em $C$-real} (or
briefly {\em real}) if $Cf=f$. Set
   \begin{align*}
\rec f=\frac{1}{2}(f + Cf) \quad \text{and} \quad \imc
f= \frac{1}{2\I}(f-Cf), \quad f\in \hh.
   \end{align*}
Then clearly for every $f\in \hh$,
   \begin{align} \label{cartC1}
\text{$\rec f$ and $\imc f$ are $C$-real, and $f =
\rec f + \I \cdot \imc f$.}
   \end{align}
Hence $\is{\rec f}{\imc f} = \is{C(\rec f)}{C(\imc
f)}=\is{\imc f}{\rec f}$ for all $f \in \hh$, and thus
$\is{\rec f}{\imc f} \in \rbb$ for all $f \in \hh$,
which gives
   \begin{align*}
\|f\|^2 = \|\rec f\|^2 + \|\imc f\|^2, \quad f \in
\hh.
   \end{align*}
Recall that if $A$ is a symmetric operator in $\hh$
such that $CA \subseteq AC$ (i.e., $C(\dz{A})
\subseteq \dz{A}$ and $CA f= ACf$ for all $f \in
\dz{A}$), then $CA^* \subseteq A^*C$, i.e.,
   \begin{align} \label{bs3}
\text{$C(\dz{A^*}) \subseteq \dz{A^*}$ and $CA^* f=
A^*Cf$ for $f \in \dz{A^*}$.}
   \end{align}
For much of the rest of the paper we will be
considering the following situation.
   \begin{align}   \label{ass}
   \begin{minipage}{72ex} Let $A$ be a closed symmetric
operator in a complex Hilbert space $\hh$ such that
$\ker(A)=\{0\}$. Suppose that there exists a complex
conjugation $C$ on $\hh$ such that $A$ is {\em
$C$-real} (or briefly {\em real}), i.e., $CA \subseteq
AC$.
   \end{minipage}
   \end{align}
   The next two lemmata are of technical importance.
   \begin{lem} \label{bs1}
Suppose that \eqref{ass} holds and $\hh \neq \{0\}$.
Then there exists $f \in \dz{A}$ such that either
$\is{f}{Af} > 0$ or $\is{f}{Af} < 0$. In the former
case, there exists $h\in \dz{A}$ such that $Ch=h$ and
$\is{h}{Ah} > 0$. In the latter case, there exists
$h\in \dz{A}$ such that $Ch=h$ and $\is{h}{Ah} < 0$.
   \end{lem}
   \begin{proof}
Since the possibility that $\is{f}{Af} = 0$ for all $f
\in \dz{A}$ is excluded by the fact that $\hh \neq
\{0\}$ and $\ker(A)=\{0\}$, and $\is{f}{Af} \in \rbb$
for all $f \in \dz{A}$, it remains to prove the last
two statements of the conclusion. By symmetry, it
suffices to consider the case when $\is{f}{Af} > 0$
for some $f \in \dz{A}$. Since $A$ is $C$-real, we
deduce that $u:=\rec f \in \dz{A}$, $v:= \imc f \in
\dz{A}$, $Cu=u$, $Cv=v$ and
   \begin{align*}
\is{u}{Av} = \is{CAv}{Cu} = \is{ACv}{u} = \is{Av}{u} =
\overline{\is{u}{Av}},
   \end{align*}
which together with $A\subseteq A^*$ implies that
   \begin{align*}
0 < \is{u + \I v}{A(u + \I v)} = \is{u}{Au} + 2
\mathrm{Re} (\I \is{u}{Av}) + \is{v}{Av} = \is{u}{Au}
+ \is{v}{Av}.
   \end{align*}
Therefore either $\is{u}{Au}>0$ or $\is{v}{Av} > 0$,
which completes the proof.
   \end{proof}
   \begin{lem} \label{BS2.5}
If $T$ is an operator in $\hh$ and $C$ is a complex
conjugation on $\hh$ such that $CT \subseteq TC$ and
$\ker(T) \neq \{0\}$, then there exists $f \in \ker(T)
\setminus \{0\}$ such that $Cf=f$.
   \end{lem}
   \begin{proof}
Take $f\in \ker(T) \setminus \{0\}$. Since $CT
\subseteq TC$ implies $C(\ker(T)) = \ker(T)$, we get
$\rec f, \imc f \in \ker{T}$, which together with
\eqref{cartC1} gives that either $\rec f \neq 0$ or
$\imc f \neq 0$. This completes the proof.
   \end{proof}
\subsection{}
Now, we concentrate on the class of vectors satisfying
the assumptions of the independence assertion of
Theorem \ref{BSth}. Given a symmetric operator $A$ in
$\hh$, a complex conjugation $C$ on $\hh$ and a vector
$\varphi$ in $\ker(A^*)$, we write
   \begin{align*}
\sac = \Big\{\psi \in \dz{A^*} \setminus \dz{A}\colon
C\psi=\psi, \is{\psi}{A^*\psi}=0,
\is{\varphi}{A^*\psi} =1\Big\}.
   \end{align*}
Clearly, the equality $\is{\varphi}{A^*\psi} =1$
implies that
   \begin{align} \label{bsr1}
\sac = \Big\{\psi \in \dz{A^*}\colon C\psi=\psi,
\is{\psi}{A^*\psi}=0, \is{\varphi}{A^*\psi} =1\Big\}.
   \end{align}
   \begin{rem}
Note that if $\varphi\in \ker(A^*)$ and $\psi \in
\dz{A^*}$, then
   \begin{align*}
\is{\alpha\varphi + \beta\psi}{A^* (\gamma\varphi +
\delta \psi)} = \left\langle \left[\begin{matrix}
\alpha \\ \beta \end{matrix} \right],
\left[\begin{matrix} 0 &
\is{\varphi}{A^*\psi} \\
0 & \is{\psi}{A^*\psi}\end{matrix} \right]
\left[\begin{matrix} \gamma
\\ \delta \end{matrix} \right]  \right\rangle, \quad
\alpha, \beta, \gamma, \delta \in \cbb.
   \end{align*}
Hence, if additionally $\varphi \neq 0$, then
$\is{\varphi}{A^*\psi} =1$ and $\is{\psi}{A^*\psi}=0$
if and only if the vectors $\varphi,\psi$ are linearly
independent and $\is{\cdot}{A^* \cdot}$ has the matrix
representation (2.7) in the basis $(\varphi, \psi)$.
   \end{rem}
   The following lemma is a modified version of what
can be found in \cite[Theorem 2.6]{sim}. For the
reader's convenience we include its proof.
   \begin{lem} \label{BS3}
Suppose that \eqref{ass} holds, $d_+(A)=1$ and
$\dim\ker(A^*) = 1$. Let $\varphi$ be a $C$-real
vector in $\ker(A^*) \setminus \{0\}$ $($cf.\ Lemma
{\em \ref{BS2.5}}$)$. Then $\sac \neq \varnothing$.
   \end{lem}
   \begin{proof}
Since $A$ is $C$-real, we infer from the von Neumann
theorem that $d_-(A) = d_+(A)=1$. This and the
equality
   \begin{align*}
\text{$\dz{A^*} = \dz{A} \dotplus \ker(A^* +\I I)
\dotplus \ker(A^* - \I I)$ \quad (direct sum),}
   \end{align*}
which is true for arbitrary closed symmetric operators
(cf.\ \cite[Lemma, p.\ 138]{R-S}), imply that $\dim
\big(\dz{A^*}/\dz{A}\big)=2$. Since $A \subseteq A^*$
and $\ker(A)=\{0\}$, we get $\dz{A}\cap \ker{A^*}
=\{0\}$. Hence $\dim \big[(\dz{A} \dotplus
\ker(A^*))/\dz{A}\big]=1$ (because $\dim
\ker(A^*)=1$), and thus there exists $\eta \in
\dz{A^*} \setminus (\dz{A} \dotplus \ker(A^*))$.
Since, by \eqref{bs3}, the vectors $\rec \eta$ and
$\imc \eta$ are in $\dz{A^*}$, we deduce from
\eqref{cartC1} that either $\rec \eta \notin \dz{A}
\dotplus \ker(A^*)$ or $\imc \eta\notin \dz{A}
\dotplus \ker(A^*)$. Therefore, we can assume without
loss of generality that $C \eta = \eta$. Putting all
these together, we get
   \begin{align} \label{dza*}
\dz{A^*} = \dz{A} \dotplus \ker(A^*) \dotplus \cbb
\cdot \eta.
   \end{align}
Now we show that $\is{\varphi}{A^*\eta} \neq 0$.
Suppose that contrary to our claim
$\is{\varphi}{A^*\eta} = 0$. Define a sesquilinear
form $Q$ on $\dz{A^*}$ by $Q(f,g)=\is{f}{A^*g} -
\is{A^*f}{g}$ for $f,g \in \dz{A^*}$. Since $A$ is
symmetric, we have $Q(f,g)=0$ for $f,g \in \dz{A}$.
Note also that
   \begin{align}  \label{ea*e}
\is{A^*\eta}{\eta} = \is{C\eta}{CA^*\eta}
\overset{\eqref{bs3}}= \is{\eta}{A^*C\eta} =
\is{\eta}{A^*\eta}.
   \end{align}
Thus $Q(\eta,\eta)=0$. Using \eqref{dza*} and
$\is{\varphi}{A^*\eta} = 0$, it is now easily seen
that $Q\equiv 0$, which means that $A^*$ is symmetric.
This and $A=\bar A$ imply that $A$ is selfadjoint,
which contradicts $d_+(A)=1$, and finally shows that
$\is{\varphi}{A^*\eta} \neq 0$. Since $\varphi$ and
$\eta$ are $C$-real, we infer from \eqref{bs3} that
$\is{\varphi}{A^*\eta} \in \rbb$. Therefore, we can
assume without loss of generality that
$\is{\varphi}{A^*\eta}=1$. Now, by setting $\psi=\eta-
\is{\eta}{A^*\eta} \varphi$, we infer from
\eqref{bsr1} and \eqref{ea*e} that $\psi \in \sac$
(our particular choice of $\psi$ guarantees that $\psi
\notin \dz{A} \dotplus \ker(A^*)$).
   \end{proof}
The next lemma is a main ingredient of the proof of
Proposition \ref{BS5}.
   \begin{lem} \label{BS4}
Suppose that \eqref{ass} holds, $d_+(A)=1$ and
$\dim\ker(A^*) = 1$. Let $\varphi$ and $h$ be $C$-real
vectors such that $\varphi \in \ker(A^*) \setminus
\{0\}$, $h\in \dz{A}$ and $\is{h}{Ah} \neq 0$ $($cf.\
Lemmata {\em \ref{bs1}} and {\em \ref{BS2.5}}$)$, and
let $\psi \in \sac$ $($cf.\ Lemma {\em \ref{BS3}}$)$.
Set
   \begin{align} \label{2def}
\text{$\widehat{\psi}(x) = \widehat{\eta}(x) -
\is{\widehat{\eta}(x)}{A^*\widehat{\eta}(x)} \varphi$
\, with \, $\widehat{\eta}(x) = x h + \varphi + \psi$
\, for $x \in \rbb$.}
   \end{align}
Then $\big\{\widehat{\psi}(x)\colon x \in \rbb\big\}
\subseteq \sac$ and
   \begin{align} \label{bsl5}
\widehat{\psi}(x) - \widehat{\psi}(y) = (x-y)h -
(\varDelta(x) - \varDelta(y))\varphi, \quad x,y \in
\rbb,
   \end{align}
where $\varDelta(x) = \varDelta_{h,\psi}(x) := x^2
\is{h}{Ah} + 2 x \mathrm{Re} \is{\psi}{Ah}$ for $x\in
\rbb$. Moreover, for every $\vartheta \in
\rbb\setminus \{0\}$ there exist $x,y\in \rbb$ such
that
   \begin{align}  \label{bsl8}
\widehat{\psi}(x) - \widehat{\psi}(y) = (x-y)h +
\vartheta \varphi \in (\dz{A} \dotplus \ker(A^*))
\setminus \dz{A}.
   \end{align}
   \end{lem}
   \begin{proof}
Since $h, \varphi, \psi$ are $C$-real, so are
$\widehat{\eta}(x)$, $x\in \rbb$. It follows from
$\varphi \perp \ob{A}$, $\is{\varphi}{A^*\psi} =1$ and
$\is{\psi}{A^*\psi}=0$ that
   \begin{align} \label{bsl1}
\is{\widehat{\eta}(x)}{A^*\widehat{\eta}(x)} = \is{x h
+ \varphi + \psi}{xAh + A^*\psi} = \varDelta(x) +1,
\quad x \in \rbb.
   \end{align}
These two facts imply that $C\widehat{\psi}(\cdot) =
\widehat{\psi}(\cdot)$. As $\varphi \perp \ob{A}$, we
have for all $x \in \rbb$
   \begin{align} \label{bs5}
\is{\varphi}{A^*\widehat{\psi}(x)} =
\is{\varphi}{A^*\widehat{\eta}(x)} = x\is{\varphi}{Ah}
+ \is{\varphi}{A^*\psi} = 1.
   \end{align}
Using the fact that $\widehat{\eta}(x)$ is $C$-real
for all $x\in \rbb$ and arguing as in \eqref{ea*e}, we
see that $\is{\widehat{\eta}(x)}{A^*
\widehat{\eta}(x)} \in \rbb$ for all $x\in \rbb$. This
yields
   \begin{align*}
\is{\widehat{\psi}(x)}{A^* \widehat{\psi}(x)} =
\is{\widehat{\eta}(x)}{A^* \widehat{\eta}(x)} -
\overline{\is{\widehat{\eta}(x)}{A^*
\widehat{\eta}(x)}} \is{\varphi}{A^*\widehat{\eta}(x)}
\overset{\eqref{bs5}}= 0, \quad x \in \rbb.
   \end{align*}
Hence $\big\{\widehat{\psi}(x)\colon x \in \rbb\big\}
\subseteq \sac$. Noting that for every $x \in \rbb$,
   \begin{align} \label{bsl-ad}
\widehat{\psi}(x) \overset{\eqref{2def}}= xh + (1 -
\is{\widehat{\eta}(x)}{A^* \widehat{\eta}(x)}) \varphi
+ \psi \overset{\eqref{bsl1}} = x h - \varDelta(x)
\varphi + \psi,
   \end{align}
we obtain \eqref{bsl5}. The latter together with
$\is{h}{Ah} \neq 0$ and the equality
   \begin{align*}
\varDelta(x) - \varDelta(y) = (x^2-y^2)\is{h}{Ah} + 2
(x-y) \mathrm{Re} \is{\psi}{Ah}, \quad x,y \in \rbb,
   \end{align*}
imply the moreover part of the conclusion.
   \end{proof}
\subsection{}
Clearly, the description of the operator $B_\infty$
given in Theorem \ref{BSth} does not depend on the
choice of $\psi\in \sac$. However, as shown in
Proposition \ref{BS5} below, the description of the
operators $\{B_t\colon t\in \rbb\}$ is extremely
dependent on the choice of the vector $\psi\in \sac$.
   \begin{pro}\label{BS5}
Suppose that \eqref{ass} holds, $d_+(A)=1$ and
$\dim\ker(A^*) = 1$. Let $\varphi$ be a $C$-real
vector in $\ker(A^*) \setminus \{0\}$ $($cf.\ Lemma
{\em \ref{BS2.5}}$)$. Then for every $(t_1,t_2) \in
\rbb \times \rbb$ such that $t_1\neq t_2$, there exist
$(\psi_1,\psi_2) \in \sac \times \sac$ such that
   \begin{align*}
   \dtf{t_1}{\psi_1} = \dtf{t_2}{\psi_2},
   \end{align*}
where $\dtf{t}{\psi} = \dz{A} + \{\alpha (t\varphi +
\psi) \colon \alpha \in \cbb\}$ for $t\in \rbb$ and
$\psi \in \sac$.
   \end{pro}
   \begin{proof}
Take $(t_1,t_2) \in \rbb \times \rbb$ such that
$t_1\neq t_2$, and fix $\psi \in \sac$ (cf.\ Lemma
\ref{BS3}). Let $h$ be a $C$-real vector in $\dz{A}$
such that $\is{h}{Ah} \neq 0$ (cf.\ Lemma \ref{bs1}),
and let $\widehat{\psi}(\cdot)$ be as in \eqref{2def}.
Set $\vartheta = t_2 - t_1$. Then by Lemma \ref{BS4},
there exist $x,y\in \rbb$ such that
$\psi_1:=\widehat{\psi}(x)\in\sac$,
$\psi_2:=\widehat{\psi}(y) \in \sac$ and \eqref{bsl8}
holds. Since $h \in \dz{A}$, we have
   \begin{align*}
\dtf{t_1}{\psi_1} & \overset{\eqref{bsl8}}= \dz{A} +
\cbb \cdot \big(t_1\varphi + \psi_2 + (x-y)h +
\vartheta \varphi)
   \\
& \hspace{1.5ex}= \dz{A} + \cbb \cdot
\big((t_1+\vartheta)\varphi + \psi_2) =
\dtf{t_2}{\psi_2},
   \end{align*}
which completes the proof.
   \end{proof}

Calculating the vectors $\widehat{\psi}(x) - \psi$
with the help of \eqref{bsl-ad} and considering them
instead of $\widehat{\psi}(x) - \widehat{\psi}(y)$ in
the proof of Proposition \ref{BS5}, we obtain the
following result which itself implies Proposition
\ref{BS5} (note that by Lemma \ref{bs1} there is no
loss of generality in assuming that the vector $h$ in
Proposition \ref{BS6} below is $C$-real).
   \begin{pro} \label{BS6}
Suppose that the assumptions of Proposition {\em
\ref{BS5}} are satisfied and $(\psi,t)\in \sac\times
\rbb$. If there exists $h\in \dz{A}$ such that
$\is{h}{Ah}>0$ $($respectively, $\is{h}{Ah}<0$$)$,
then for every real $t^\prime > t$ $($respectively,
$t^\prime < t$$)$, there exists $\psi^\prime \in \sac$
such that $\dtf{t}{\psi} =
\dtf{t^\prime}{\psi^\prime}$.
   \end{pro}
The following proposition together with Lemma
\ref{BS4} shows that the term $\ker(A^*)$ which
appears in the formula \eqref{formulaF1} below could
not be removed without spoiling the conclusion of
Proposition \ref{F1} (in contrast to what is written
in the proof of the independence assertion of
\cite[Theorem 2.6]{sim}).
   \begin{pro}\label{F1}
Suppose that the assumptions of Proposition {\em
\ref{BS5}} are satisfied. If $\psi$ and $\psi^\prime$
are any two vectors in $\sac$, then
   \begin{align} \label{formulaF1}
\psi^\prime - \psi \in \dz{A} \dotplus \ker(A^*)
=\dz{A} \dotplus \cbb \cdot \varphi.
   \end{align}
   \end{pro}
   \begin{proof}
First we note that $\psi \notin \dz{A} \dotplus \cbb
\cdot \varphi$. Indeed, otherwise $\psi = f + \gamma
\cdot \varphi$ for some $f \in \dz{A}$ and $\gamma \in
\cbb$, which implies
   \begin{align*}
1 \overset{\eqref{bsr1}}= \is{\varphi}{A^*\psi} =
\is{\varphi}{Af} = \is{A^*\varphi}{f} = 0,
   \end{align*}
a contradiction. Since $\dim
\big(\dz{A^*}/\dz{A}\big)=2$ (see the proof of Lemma
\ref{BS3}), we deduce that $\dz{A^*} = \dz{A} \dotplus
\cbb \cdot \varphi \dotplus \cbb \cdot \psi$. Hence,
there exist $h \in \dz{A}$ and $\alpha,\beta \in \cbb$
such that $\psi^\prime - \psi = h + \alpha\cdot
\varphi + \beta \cdot \psi$, which yields
   \begin{align*}
0 \overset{\eqref{bsr1}}= \is{A^*\psi^\prime}{\varphi}
- \is{A^*\psi}{\varphi} = \is{A^*(\psi^\prime -
\psi)}{\varphi} = \is{Ah}{\varphi} + \bar\beta
\is{A^*\psi}{\varphi} \overset{\eqref{bsr1}}=
\bar\beta.
   \end{align*}
Thus $\psi^\prime - \psi = h + \alpha\cdot \varphi \in
\dz{A} \dotplus \cbb \cdot \varphi$, which together
with the equality $\ker(A^*)=\cbb \cdot \varphi$
completes the proof.
   \end{proof}
The question of when two vectors $\psi, \psi^\prime
\in \sac$ represent the same operators $\{B_t\colon
t\in \rbb\}$ in the sense that $\dtf{t}{\psi} =
\dtf{t}{\psi^\prime}$ for all $t \in \rbb$ has a
simple answer.
   \begin{pro} \label{F2}
Suppose that the assumptions of Proposition {\em
\ref{BS5}} are satisfied. If $\psi, \psi^\prime \in
\sac$, then the following conditions are
equivalent{\em :}
   \begin{enumerate}
   \item[(i)] $\dtf{t}{\psi} = \dtf{t}{\psi^\prime}$
for all $t\in \rbb$,
   \item[(ii)] there exist $\alpha, \beta \in \cbb$ such
that $|\alpha|^2+|\beta|^2 >0$ and $\alpha\psi + \beta
\psi^\prime \in \dz{A}$,
   \item[(iii)] $\psi^\prime -
\psi \in \dz{A}$.
   \end{enumerate}
   \end{pro}
   \begin{proof}
(i)$\Rightarrow$(ii) Since $\psi^\prime \in
\dtf{0}{\psi} = \dtf{0}{\psi^\prime}$, we see that
$\psi^\prime = h + \alpha \psi$ for some $h\in \dz{A}$
and $\alpha \in \cbb$.

(ii)$\Rightarrow$(iii) Since $h:=\alpha\psi + \beta
\psi^\prime \in \dz{A}$, $\varphi \in \ker(A^*)$ and
$\psi, \psi^\prime \in \sac$, we get
   \begin{align*}
0 = \is{\varphi}{A^*(\alpha\psi + \beta \psi^\prime)}
= \alpha\is{\varphi}{A^*\psi} + \beta
\is{\varphi}{A^*\psi^\prime} \overset{\eqref{bsr1}}=
\alpha + \beta.
   \end{align*}
Hence, by the inequality $|\alpha|^2+|\beta|^2 >0$, we
have $\beta\neq 0$ and $\psi^\prime - \psi =
\beta^{-1}h \in \dz{A}$.

(iii)$\Rightarrow$(i) Obvious.
   \end{proof}
   \begin{rem}
In view of the above discussion it is natural to ask
whether the following implication is valid (still
under the assumptions of Proposition \ref{BS5}):
   \begin{align} \label{iMpl}
\psi \in \sac, \psi^\prime \in \dz{A^*},
C\psi^\prime=\psi^\prime, \psi^\prime - \psi \in
\dz{A} \implies \psi^\prime \in \sac.
   \end{align}
We show that the answer is in the negative (note,
however, that $\is{\varphi}{A^* \psi^\prime}=1$).
Suppose that, contrary to our claim, the implication
\eqref{iMpl} is valid. Take $\psi \in \sac$ and a
$C$-real vector $h \in \dz{A}$. Then, by \eqref{iMpl}
applied to $\psi^\prime = t h + \psi$, we obtain
   \begin{align*}
0=\is{t h + \psi}{A^*(t h + \psi)}=t^2\is{h}{Ah} + 2 t
\mathrm{Re} \is{h}{A^*\psi}, \quad t \in \rbb.
   \end{align*}
Hence $\is{h}{Ah}=0$ for all $C$-real vectors $h \in
\dz{A}$. This contradicts Lemma \ref{bs1}.
   \end{rem}
   \textbf{Acknowledgement}. A substantial part of
this paper was written while the first and the third
authors visited Kyungpook National University during
the spring and the autumn of 2011; they wish to thank
the faculty and the administration of this unit for
their warm hospitality.
   \bibliographystyle{amsalpha}
   
   \end{document}